\documentclass[10pt]{amsart}
\usepackage[dvipdfmx]{graphicx}
\usepackage{amsmath,amscd,amsthm,amsxtra,yhmath}
\usepackage{epsfig,graphics,color,colortbl}
\usepackage{amssymb,latexsym}
\usepackage{mathrsfs,upgreek}
\usepackage[poly,all]{xy}
\usepackage{hyperref}
\usepackage{tikz-cd}
\usepackage{ytableau}

\newtheorem{thm}{\bf Theorem}[section]
\newtheorem{df}[thm]{\bf Definition}
\newtheorem{prop}[thm]{\bf Proposition}
\newtheorem{cor}[thm]{\bf Corollary}
\newtheorem{lem}[thm]{\bf Lemma}
\newtheorem{rem}[thm]{\bf Remark}
\newtheorem{ex}[thm]{\bf Example}

\numberwithin{equation}{section}

\newcommand{\mc}{\mathcal}
\newcommand{\mf}{\mathfrak}

\newcommand{\pf}{\noindent{\bfseries Proof. }}

\newcommand{\U}{{\mc U}}
\newcommand{\V}{\mc{V}}
\newcommand{\W}{\mc{W}}
\newcommand{\cP}{\mathscr{P}}
\newcommand{\cO}{\mc{O}}
\newcommand{\I}{\mathbb{I}}
\newcommand{\be}{{\bf e}}

\newcommand{\Z}{\mathbb{Z}}
\newcommand{\Q}{\mathbb{Q}}

\newcommand{\e}{\epsilon}
\newcommand{\de}{\delta}
\newcommand{\te}{\widetilde{e}}
\newcommand{\tf}{\widetilde{f}}
\newcommand{\td}{\widetilde}

\newcommand{\la}{\lambda}

\newcommand{\si}{(-1)^{\e_i}}

\newcommand{\ot}{\otimes}

\begin{document}
\title[]
{$R$ matrix for generalized quantum group of type $A$}

\author{JAE-HOON KWON}

\address{Department of Mathematical Sciences and Research Institute of Mathematics, Seoul National University, Seoul 08826, Korea}
\email{jaehoonkw@snu.ac.kr}

\author{JEONGWOO YU}

\address{Department of Mathematical Sciences, Seoul National University, Seoul 08826, Korea}
\email{ycw453@snu.ac.kr}


\thanks{J.-H.K. was supported by the National Research Foundation of Korea(NRF) grant funded by the Korea government(MSIT) (No. 2019R1A2C108483311).}

\begin{abstract}
The generalized quantum group $\U(\e)$ of type $A$ is an affine analogue of quantum group associated to a general linear Lie superalgebra $\mf{gl}_{M|N}$.
We prove that there exists a unique $R$ matrix on tensor product of fundamental type representations of $\U(\e)$ for arbitrary parameter sequence $\e$ corresponding to a non-conjugate Borel subalgebra of $\mf{gl}_{M|N}$.
We give an explicit description of its spectral decomposition, and then as an application,  construct a family of finite-dimensional irreducible $\U(\e)$-modules which have subspaces isomorphic to the Kirillov-Reshetikhin modules of usual affine type $A_{M-1}^{(1)}$ or $A_{N-1}^{(1)}$.

\end{abstract}

\maketitle
\setcounter{tocdepth}{1}
\tableofcontents

\section{Introduction}

A generalized quantum group $\U(\e)$ associated to $\epsilon=(\e_1,\ldots,\e_n)$ with $\e_i\in\{0,1\}$ is a Hopf algebra introduced in \cite{KOS}, which appears in the study of solutions to the tetrahedron equation or the three-dimensional Yang-Baxter equation.

The generalized quantum group $\U(\e)$ of type $A$ is equal to the usual quantum affine algebra of type $A_{n-1}^{(1)}$, when $\e$ is homogeneous, that is, $\epsilon_i=\epsilon_j$ for all $i\neq j$. 
But it becomes a more interesting object when $\e$ is non-homogeneous, which is closely related to the quantized enveloping algebra associated to an affine Lie superalgebra \cite{Ya99}, or which can be viewed as an affine analogue of the quantized enveloping algebra of the general linear Lie superalgebra $\mf{gl}_{M|N}$ \cite{Ya94}, where $M$ and $N$ are the numbers of $0$ and $1$ in $\e$, respectively. We remark that the subalgebra $\ring\U(\e)$ of $\U(\e)$ associated to the Lie superalgebra $\mf{gl}_{M|N}$ was also introduced in \cite{CWZ} independently, as symmetries appearing in the study of wave functions of quantum mechanical systems \cite{Za}.

When the parameter $\e$ is standard, that is, $\e_{M|N}=(0^M,1^N)$, it is shown in \cite{KOS} that there exists a unique $R$ matrix on the tensor product of finite-dimensional $\U(\e_{M|N})$-modules $\W_{s,\e}(x)$, which correspond to fundamental representations of type $A_{N-1}^{(1)}$ with spectral parameter $x$ when $N\geq 3$.
Indeed, the $R$ matrix is obtained by reducing the solution of the tetrahedron equation, and the uniqueness follows from the irreducibility of tensor product $\W_{l,\e}(x)\ot \W_{m,\e}(y)$ for generic $x$ and $y$. An explicit spectral decomposition of the associated $R$ matrix is obtained by analyzing the maximal vectors with respect to $\ring\U(\e_{M|N})$. 

By applying fusion construction using the $R$ matrix in \cite{KOS}, a family of irreducible $\U(\e_{M|N})$-modules is constructed in \cite{KO}, which are parametrized by rectangular partitions inside an $(M|N)$-hook. Moreover the existence of their crystal base is proved together with a combinatorial description of the associated crystal graphs. It can be viewed as a natural super-analogue of Kirillov-Reshetikhin modules (simply KR modules) of type $A_{\ell}^{(1)}$, which is a most important family of finite-dimensional irreducible modules of quantum affine algebras (cf.~\cite{CH,KR2}).

The results in \cite{KOS} and \cite{KO} suggests that there is a close connection between finite-dimensional representations of $\U(\e_{M|N})$ and $U_q(A_{n-1}^{(1)})$.
The purpose of this paper is to extend the results in \cite{KOS} and \cite{KO} to arbitrary parameter sequence $\e$, and find a more concrete connection between the finite-dimensional representations of $\U(\e)$ and  $U_q(A_{\ell}^{(1)})$. From a viewpoint of representations of $\mf{gl}_{M|N}$, the sequence $\e$ represents the type of Borel subalgebras of $\mf{gl}_{M|N}$, which are not conjugate to each other. It is not obvious whether the representation theory of $\U(\e)$  is the same under a different choice of permutations of $\e_{M|N}$. For example, if we change the Borel in the generalized quantum group, then the defining relations and the crystal structure associated to $\U(\e)$-modules become much different from the ones with respect to $\e_{M|N}$ as $\e$ gets far from $\e_{M|N}$ (cf.~\cite{BKK,K09}).

We first show that there exists a unique $R$ matrix on the tensor product of finite-dimensional $\U(\e)$-modules $\W_{s,\e}(x)$ of fundamental type (Theorem \ref{thm:existence of R-matrix}). 
Since the existence of $R$ matrix for arbitrary $\e$ was shown in \cite{KOS}, it suffices to prove the irreducibility of tensor product $\W_{l,\e}(x)\ot \W_{m,\e}(y)$ for generic $x$ and $y$. We use a method completely different from \cite{KOS}.
Indeed, motivated by the work \cite{CL}, we introduce a functor called {truncation}, and show that it sends any $\U(\e)$-module with polynomial weights to a $\U(\e')$-module, preserving the comultiplications in tensor product, where $\e'$ is a subsequence of $\e$. This in particular enables us to define an oriented graph structure on $\W_{l,\e}(x)\ot \W_{m,\e}(y)$ when $x=y=1$ with additional arrows other than the ones associated to $\U(\e)$. With this structure, we prove the connectedness of the crystal (Theorem \ref{thm:main theorem-1}), and hence the irreducibility for generic $x$ and $y$.

Next, we prove that the truncation functor is compatible with the $R$ matrix. 
This immediately implies that the spectral decomposition of the $R$ matrix for $\U(\e)$ is the same as that of type $A_{\ell}^{(1)}$ (Theorem \ref{thm:spectral decomp}) and hence does not depend on the choice of $\e$.
As an application, we construct a family of irreducible $\U(\e)$-modules $\W^{(r)}_{s,\e}$ which  yields the usual KR modules under truncation (Theorem \ref{thm:W_s^{(r)}}). We conjecture that $\W^{(r)}_{s,\e}$ has a crystal base as in the case of $\e=\e_{M|N}$. We expect that the compatibility of truncation with the $R$ matrix will also play a crucial role in understanding arbitrary finite-dimensional $\U(\e)$-modules in connection with those of type $A_\ell^{(1)}$.

There are other recent works on the finite-dimensional representations of quantum affine {\em superalgebra} associated to $\mf{gl}_{M|N}$ \cite{Z14,Z16,Z17}. It would be interesting to compare with these results.

The paper is organized as follows. In Section \ref{sec:GQG}, we review basic materials for a generalized quantum group and its crystal base. In Section \ref{sec:poly reps}, we present the classical Schur-Weyl duality for $\ring\U(\e)$ and then realize the irreducible polynomial representation of $\ring\U(\e)$. In Section \ref{sec:R matrix}, we prove the main theorem on the existence of the $R$ matrix. In Section \ref{sec:KR modules}, we construct KR type modules of $\U(\e)$ using the $R$ matrix.\vskip 3mm

{\bf Acknowledgement} The authors would like to thank Euiyong Park and Masato Okado for helpful discussions and thank Shin-Myung Lee for careful reading of the manuscript and comments.

\section{Generalized quantum group $\U(\e)$ of type $A$}\label{sec:GQG}

\subsection{Generalized quantum group}
We fix a positive integer $n\geq 4$.
Let $\e=(\e_1,\cdots,\e_n)$ be a sequence with $\e_i\in \{0,1\}$ for $1\leq i\leq n$. We denote by $\I$ the linearly ordered set $\{1<2<\cdots <n\}$ with $\Z_2$-grading given by $\I_0=\{\,i\,|\,\e_i=0\,\}$ and $\I_1=\{\,i\,|\,\e_i=1\,\}$.
We assume that $M$ is the number of $i$ with $\e_i=0$ and $N$ is the number of $i$ with $\e_i=1$ in $\e$. We denote by $\e_{M|N}$ the sequence when $\e_1=\dots=\e_M=0$ and $\e_{M+1}=\dots=\e_n=1$.

Let $P=\bigoplus_{i\in \I}\Z\delta_i$ be the free abelian group generated by $\de_i$ with a symmetric bilinear form $(\,\cdot\,|\,\cdot\,)$ given by $(\de_i|\de_j)=(-1)^{\e_i}\de_{ij}$ for $i,j\in \I$. Let $\{\,\de^\vee_i\,|\,i\in \I\,\}\subset P^\vee:={\rm Hom}_\Z(P,\Z)$ be the dual basis such that $\langle \de_i, \de^\vee_j \rangle =\de_{ij}$ for $i,j\in \I$.

Let $I=\{\,0,1,\ldots,n-1\,\}$ and
\begin{equation*}\label{eq:simple root}
\alpha_i=\de_i-\de_{i+1},\quad \alpha_i^\vee = \de^\vee_i-(-1)^{\e_i+\e_{i+1}}\de^\vee_{i+1} \quad (i\in I).
\end{equation*}
Throughout the paper, we understand the subscript $i\in I$ modulo $n$.
When $\e=\e_{M|N}$, the Dynkin diagram associated to the Cartan matrix $(\langle \alpha_j,\alpha^\vee_i \rangle)_{0\leq i,j\leq n}$ is  
\vskip 4mm
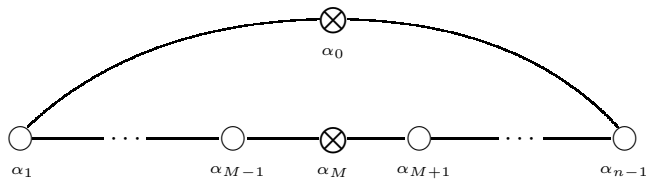
\begin{center}
\hskip -4cm \setlength{\unitlength}{0.25in}
\begin{picture}(24,4)
\put(7.95,1){\makebox(0,0)[c]{$\bigcirc$}}
\put(10.2,1){\makebox(0,0)[c]{$\cdots$}}
\put(12.4,1){\makebox(0,0)[c]{$\bigcirc$}}
\put(16.3,1){\makebox(0,0)[c]{$\bigcirc$}}
\put(18.45,1){\makebox(0,0)[c]{$\cdots$}}
\qbezier(8.1,1.25)(10.4,3.4)(14.2,3.5)
\qbezier(14.8,3.5)(18.45,3.4)(20.45,1.25)
\put(8.2,1){\line(1,0){1.5}}
\put(10.6,1){\line(1,0){1.5}}
\put(12.7,1){\line(1,0){1.5}}
\put(14.8,1){\line(1,0){1.2}}
\put(16.55,1){\line(1,0){1.4}}
\put(18.9,1){\line(1,0){1.4}}
\put(14.5,1){\makebox(0,0)[c]{$\bigotimes$}}
\put(14.5,3.5){\makebox(0,0)[c]{$\bigotimes$}}
\put(20.6,1){\makebox(0,0)[c]{$\bigcirc$}}
\put(14.5,2.8){\makebox(0,0)[c]{\tiny $\alpha_{0}$}}
\put(14.5,0.3){\makebox(0,0)[c]{\tiny $\alpha_{M}$}}
\put(8,0.3){\makebox(0,0)[c]{\tiny $\alpha_{1}$}}
\put(20.6,0.3){\makebox(0,0)[c]{\tiny $\alpha_{n-1}$}}
\put(12.5,0.3){\makebox(0,0)[c]{\tiny $\alpha_{M-1}$}}
\put(16.4,0.3){\makebox(0,0)[c]{\tiny $\alpha_{M+1}$}}
\end{picture}
\end{center}%
\noindent where $\bigotimes$ denotes an isotropic simple root.

Let $q$ be an indeterminate. We put
$I_{\rm even}=\{\,i\in I\,|\,(\alpha_i|\alpha_i)=\pm 2\,\}$ and 
$I_{\rm odd}=\{\,i\in I\,|\,(\alpha_i|\alpha_i)=0\,\}$,
%
and set 
\begin{equation*}
q_i=\si q^{\si}=
\begin{cases}
q & \text{if $\e_i=0$},\\
-q^{-1} & \text{if $\e_i=1$},\\
\end{cases} \quad (i\in I).
\end{equation*}
\begin{df}\label{def:U(e)}
{\rm
We define ${\U}(\e)$ to be the associative $\Q(q)$-algebra with $1$ 
generated by $q^{h}, e_i, f_i$ for $h\in P^\vee$ and $i\in I$ 
satisfying
{\allowdisplaybreaks
\begin{align}
& q^0=1, \quad q^{h +h'}=q^{h}q^{h'} \hskip 2.5cm  (h, h' \in P^{\vee}),\label{eq:Weyl-rel-1} \\ 
& \omega_je_i\omega_j^{-1}=q_{j}^{\langle\alpha_i,\delta_j^\vee\rangle}e_i,\quad \omega_jf_i\omega_j^{-1}=q_{j}^{-\langle\alpha_i,\delta_j^\vee\rangle}f_i, \label{eq:Weyl-rel-2} \\ 
&  e_if_j - f_je_i =\delta_{ij}\frac{\omega_i\omega_{i+1}^{-1} - \omega_i^{-1}\omega_{i+1}}{q-q^{-1}},\label{eq:Weyl-rel-3}\\
& e_i^2= f_i^2 =0 \hskip 4.5cm (i\in I_{\rm odd}),\label{eq:Weyl-rel-4}
\end{align}
where $\omega_j=q^{\si\de^\vee_j}$ ($j\in \I$), and the Serre-type relations
\begin{equation}\label{eq:Serre-rel-1}
\begin{split}
&\ \, e_i e_j -  e_j e_i = f_i f_j -  f_j f_i =0,
 \hskip 1cm \text{($|i-j|>1$)},\\ 
&
\begin{array}{ll}
e_i^2 e_j- (-1)^{\e_i}[2] e_i e_j e_i + e_j e_i^2= 0,\\ 
f_i^2 f_j- (-1)^{\e_i}[2] f_i f_j f_i+f_j f_i^2= 0,
\end{array}
\ \hskip 1cm\text{($i\in I_{\rm even}$ and $|i-j|=1$)}, 
\end{split}
\end{equation}
and
\begin{equation}\label{eq:Serre-rel-2}
\begin{array}{ll}
  e_{i}e_{i-1}e_{i}e_{i+1}  
- e_{i}e_{i+1}e_{i}e_{i-1} 
+ e_{i+1}e_{i}e_{i-1}e_{i} \\  
\hskip 2cm - e_{i-1}e_{i}e_{i+1}e_{i} 
+ (-1)^{\e_i}[2]e_{i}e_{i-1}e_{i+1}e_{i} =0, \\ 
  f_{i}f_{i-1}f_{i}f_{i+1}  
- f_{i}f_{i+1}f_{i}f_{i-1} 
+ f_{i+1}f_{i}f_{i-1}f_{i}  \\  
\hskip 2cm - f_{i-1}f_{i}f_{i+1}f_{i} 
+ (-1)^{\e_i}[2]f_{i}f_{i-1}f_{i+1}f_{i} =0,
\end{array}\quad \text{($i\in I_{\rm odd}$)}.
\end{equation}}
We call $\U(\e)$ the {\em generalized quantum group of affine type $A$ associated to $\e$} (see \cite{KOS}).
}
\end{df}
Put $k_i=\omega_i\omega_{i+1}^{-1}$ for $i\in I$. 
Then we have for $i,j\in I$
\begin{equation*}
k_ie_jk_i^{-1}=D_{ij}e_j,\quad 
k_if_jk_i^{-1}=D_{ij}^{-1}f_j,\quad
e_if_j - f_je_i =\delta_{ij}\frac{k_i - k_i^{-1}}{q-q^{-1}},
\end{equation*}
where $D_{ij}=q_{i}^{\langle\alpha_j,\delta_i^\vee\rangle}
q_{i+1}^{-\langle\alpha_j,\delta_{i+1}^\vee\rangle}$.
There is a Hopf algebra structure on $\U(\e)$, where the comultiplication $\Delta$, the antipode $S$, and the couint $\varepsilon$ are given by
\begin{equation}\label{eq:comult-1}
\begin{split}
\Delta(q^h)&=q^h\otimes q^h, \\ 
\Delta(e_i)&= e_i\ot 1 + k_i^{-1}\ot e_i , \\  
\Delta(f_i)&= f_i\ot k_i + 1\ot f_i, \\  
\end{split}
\end{equation}
\begin{gather*}
S(q^h)=q^{-h}, \ \ S(e_i)=-e_ik_i^{-1}, \ \  S(f_i)=-k_if_i,\\
\varepsilon(q^h)=1,\ \ \varepsilon(e_i)=\varepsilon(f_i)=0,
\end{gather*}
for $h\in P^\vee$ and  $i\in I$.
Let $\eta$ be the anti-automorphism on $\U(\e)$ defined by
\begin{equation*}\label{eq:antiauto-1}
\eta(q^h)=q^h,\quad \eta(e_i)=q_if_ik^{-1}_i,\quad \eta(f_i)=q^{-1}_ik_ie_i,
\end{equation*}
for $h \in  P^{\vee}$ and $i\in I$. It satisfies $\eta^2=id$ and 
\begin{equation*}\label{eq:antiauto-2}
 \Delta \circ \eta= (\eta \ot \eta)\circ \Delta.
\end{equation*}

We have an isomorphism between $\U(\e)$ and $\U(\td\e)$ where $\td\e$ is obtained from $\e$ by permutation of $\e_i$'s, which is not an isomorphism of Hopf algebras \cite[Theorem 2.7]{Ma} (cf. \cite[37.1]{Lu93}).

\begin{thm}\label{thm:reflection}
For $1\leq i\leq n-1$, let $\td\e=(\td\e_1,\dots,\td\e_n)$ be the sequence given by exchanging $\e_i$ and $\e_{i+1}$ in $\e$. Then there exists an isomorphism of algebras $\tau_i : \U(\e) \longrightarrow \U(\td\e)$ given by 
\begin{equation*}
\begin{split}
&\tau_i(k_i)=k_i^{-1},\quad\tau_i(e_i)=-f_ik_i,\quad \tau_i(f_i)= -k_i^{-1}e_i,\\
&\tau_i(k_j)=k_ik_j,\quad\tau_i(e_j)=[e_i,e_j]_{D_{ij}},\quad \tau_i(f_j)= [f_j,f_i]_{D_{ij}^{-1}} \quad (|i-j|=1),\\
&\tau_i(k_j)=k_j,\quad\tau_i(e_j)=e_j,\quad \tau_i(f_j)= f_j \quad (|i-j|>1),
\end{split}
\end{equation*}
where the inverse map is given by 
\begin{equation*}
\begin{split}
&\tau_i^{-1}(k_i)=k_i^{-1},\quad\tau_i^{-1}(e_i)=-k_i^{-1}f_i,\quad \tau_i^{-1}(f_i)= -e_ik_i,\\
&\tau_i^{-1}(k_j)=k_ik_j,\quad\tau_i^{-1}(e_j)=[e_j,e_i]_{D_{ij}},\quad \tau_i^{-1}(f_j)= [f_i,f_j]_{D_{ij}^{-1}} \quad (|i-j|=1),\\
&\tau_i^{-1}(k_j)=k_j,\quad\tau_i^{-1}(e_j)=e_j,\quad \tau_i^{-1}(f_j)= f_j \quad (|i-j|>1).
\end{split}
\end{equation*}
\end{thm}
\qed

\subsection{Crystal base of $\U(\e)$-modules}

For a $\U(\e)$-module $V$ and $\mu=\sum_{i}\mu_i\delta_i\in P$, let 
\begin{equation*}
V_\mu 
= \{\,u\in V\,|\,\omega_i u= q_i^{\mu_i} u \ \ (i\in \I) \,\}
\end{equation*}
be the $\mu$-weight space of $V$. For a non-zero vector $u\in V_\mu$, we denote by ${\rm wt}(u)=\mu$ the weight of $u$. 
Let $P_{\geq0}=\sum_{i\in \I}\Z_{\geq 0}\delta _i$ and 
let $\cO_{\geq0}$ be the category of $\U(\e)$-modules with objects $V$ such that
\begin{equation}\label{eq:polynomial weight}
V=\bigoplus_{\mu\in P_{\ge 0}}V_\mu \quad \text{with}\  \dim V_\mu < \infty.
\end{equation}
which is closed under taking submodules, quotients and tensor products.

\begin{rem}{\rm
There is another comultiplication on $\U(\e)$ given by
\begin{equation}\label{eq:comult-2}
\begin{split}
\Delta_+(q^h)&=q^h\otimes q^h, \\ 
\Delta_+(e_i)&= 1\ot e_i+  e_i\ot k_i, \\  
\Delta_+(f_i)&= k_i^{-1} \ot f_i + f_i\ot 1,
\end{split}
\end{equation}
(while $\Delta_+^{\rm op}$ is used in \cite{KO}).
Let $\ot$ and $\ot_+$ denote the tensor product with respect to $\Delta$ and $\Delta_+$, respectively. 
For $\U(\e)$-modules $M$ and $N$, we have a $\U(\e)$-linear isomorphism $\psi : M\ot N \longrightarrow M\ot_+ N$ given by
\begin{equation}\label{eq:isomorphism for two comult}
\psi(u\ot v) = \left(\prod_{i\in \I}q_i^{\mu_i\nu_i}\right)\, u\ot v,
\end{equation}
for $u\in M_\mu$ and $v\in N_\nu$ with $\mu=\sum_i\mu_i\de_i$ and $\nu=\sum_i\nu_i\de_i$.
}
\end{rem}

Let us recall the notion of crystal base for $V\in \cO_{\ge 0}$ \cite{KO} (cf.~\cite{BKK}). 
The Kashiwara operators $\tilde{e}_i$ and $\tilde{f}_i$ on $V$ for $i\in I$ are defined as follows. Suppose that $u\in V_\mu$ is given.

{\em Case 1}. Suppose that $i \in I_{\rm odd}$ and $(\e_i,\e_{i+1})=(0,1)$. We define
\begin{equation*}
\tilde{e}_i u =\eta(f_i) u =q_i^{-1}k_ie_i u,\quad \tilde{f}_i u=f_i u.
\end{equation*}

{\em Case 2.} Suppose that $i \in I_{\rm odd}$ and $(\e_i,\e_{i+1})=(1,0)$. We define
\begin{equation*}
\tilde{e}_i u=e_i u,\quad\tilde{f}_i u=\eta(e_i) u=q_if_ik^{-1}_iu.
\end{equation*}

{\em Case 3.} Suppose that $i \in I_{\rm even}$ and $(\e_i,\e_{i+1})=(0,0)$.
Let $\zeta: U_q(\mf{sl}_2) \rightarrow \U(\e)_i$ be the $\Q(q)$-algebra isomorphism given by
$\zeta(e)=e_i$, $\zeta(f)=f_i$ and $\zeta(k)=k_i$, where $U_q(\mf{sl}_2)=\langle e,f,k^{\pm 1} \rangle$ is the usual quantum group for $\mf{sl}_2$ with relation $kek^{-1}=q^2 e$, $kfk^{-1}=q^{-2} f$, $ef-fe=\frac{k-k^{-1}}{q-q^{-1}}$.
The induced comultiplication $\Delta^\zeta :=(\zeta^{-1} \otimes \zeta^{-1})\circ \Delta \circ \zeta$ on $U_q(\mf{sl_2})$ is 
\begin{equation*}
\begin{split}
\Delta^\zeta(k^{\pm1})&=k^{\pm1} \otimes k^{\pm1},\\
\Delta^\zeta(e)&=k^{-1} \otimes e+e\otimes1,\\
\Delta^\zeta(f)&=1 \otimes f + f \otimes k.
\end{split}
\end{equation*}
So we define $\tilde{e}_i$ and $\tilde{f}_i$ on $V$ to be the usual Kashiwara operators on the lower crystal base of $U_q(\mf{sl_2})$-module induced from $\zeta$. 
In other words, if  $u=\sum_{k \geq 0} f_i^{(k)}u_k$, where $f_i^{(k)}=f_i^k/[k]!$ and $e_iu_k=0$ for $k\geq 0$, then we define
\begin{equation*}
\tilde{e}_iu=\sum_{k\geq1}f_i^{(k-1)}u_k,\quad \tilde{f}_iu=\sum_{k\geq0}f_i^{(k+1)}u_k.
\end{equation*}

{\em Case 4.} Suppose that $i \in I_{even}$ and $(\e_i,\e_{i+1})=(1,1)$. 
Let $\xi: U_q(\mf{sl_2}) \rightarrow \U(\e)_i$ be the  $\Q(q)$-algebra homomorphism given by
$\xi(e)=-e_i$, $\xi(f)=f_i$ and $\xi(k)=k^{-1}_i$. 
Then the induced comultiplication $\Delta^\xi$ on $U_q(\mf{sl_2})$ is 
\begin{equation*}
\begin{split}
\Delta^\xi(k^{\pm1})&=k^{\pm1} \otimes k^{\pm1},\\
\Delta^\xi(e)&=k \otimes e+e\otimes1,\\
\Delta^\xi(f)&=1 \otimes f + f \otimes k^{-1}.
\end{split}
\end{equation*}
So we define $\tilde{e}_i$ and $\tilde{f}_i$ on $V$ to be the Kashiwara operators on the upper crystal base of $U_q(\mf{sl}_2)$-module induced from $\xi$. In other words, if $u=\sum_{k \geq 0} f_i^{(k)}u_k$, where  $e_iu_k=0$ for $k\geq 0$ and $l_k=\langle {\rm wt}(u_k),\alpha_i^\vee\rangle$, then we define
\begin{equation*}
\tilde{e}_iu=\sum_{k\geq1}q^{-l_k+2k-1}f_i^{(k-1)}u_k,\quad \tilde{f}_iu=\sum_{k\geq 0}q^{l_k-2k-1}f_i^{(k+1)}u_k.
\end{equation*}

Let $A_0$ be the subring of $\Q(q)$ consisting of $f(q)/g(q)$ with $f(q), g(q)\in \Q[q]$ and $g(0)\neq 0$.
\begin{df}\label{def:cry}{\rm
Let $V \in \cO_{\ge 0}$ be given. A pair $(L,B)$ is a crystal base of $V$ it if satisfies the following conditions:
\begin{itemize}
\item[(1)] $L$ is an $A_0$-lattice of $V$ and $L=\bigoplus_{\mu \in P_{\geq 0}}L_\mu$, where $L_\mu=L \cap V_\mu$,

\item[(2)] $B$ is a signed basis of $L/qL$, that is $B=\mathbf{B}\cup -\mathbf{B}$ where $\mathbf{B}$ is a $\Q$-basis of $L/qL$,

\item[(3)] $B= \bigsqcup_{\mu \in P_{\geq 0}} B_\mu$ where $B_\mu \subset (L/qL)_\mu$,

\item[(4)] $\tilde{e}_i L \subset L, \tilde{f}_i L \subset L$ and $\tilde{e}_iB \subset B \cup \{0\}, \tilde{f}_iB \subset B \cup \{0\}$ for $i \in I$,

\item[(5)] $\tilde{f}_ib=b'$ if and only if $\tilde{e}_ib'=\pm b$ for $i\in I$ and $b, b' \in B$.
\end{itemize}}
\end{df}

Let us call $B/\{\pm 1\}$ a crystal of $V$, which is an $I$-colored oriented graph.
We have a tensor product rule for crystals (see \cite{BKK} and \cite[Proposition 3.4]{KO}).

\begin{prop}\label{prop:tensor product rule}
Let $V_1, V_2\in \cO_{\geq 0}$ be given. Suppose that $(L_i,B_i)$ is a crystal base of $V_i$ for $i=1,2$. 
Then $(L_1\otimes L_2, B_1\otimes B_2)$ is a crystal base of $V_1\otimes V_2$, where $B_1\otimes B_2\subset (L_1/qL_1)\otimes (L_2/qL_2)=L_1\otimes L_2/qL_1\otimes L_2$. 
Moreover, $\te_i$ and $\tf_i$ act on $B_1\otimes B_2$ as follows:
\begin{itemize}
\item[(1)] 
if $i\in I_{\rm odd}$ and $(\e_i,\e_{i+1})=(0,1)$, then
{\allowdisplaybreaks
\begin{equation}\label{eq:tensor product rule for odd +}
\begin{split}
\te_i(b_1\otimes b_2)=&
\begin{cases}
b_1\otimes \te_i  b_2, & \text{if }\langle {\rm wt}(b_2),\alpha^\vee_i\rangle>0, \\ 
\te_i b_1\otimes  b_2, & \text{if }\langle {\rm wt}(b_2),\alpha^\vee_i\rangle=0,
\end{cases}
\\
\tf_i(b_1\otimes b_2)=&
\begin{cases}
 b_1\otimes \tf_ib_2, & \text{if }\langle {\rm wt}(b_2),\alpha^\vee_i\rangle >0, \\ 
\tf_ib_1\otimes  b_2, & \text{if }\langle {\rm wt}(b_2),\alpha^\vee_i\rangle=0,
\end{cases}
\end{split}
\end{equation}}
\item[(2)] 
if $i\in I_{\rm odd}$ and $(\e_i,\e_{i+1})=(1,0)$, then
{\allowdisplaybreaks
\begin{equation}\label{eq:tensor product rule for odd -}
\begin{split}
\te_i(b_1\otimes b_2)=&
\begin{cases}
 b_1\otimes \te_i b_2, & \text{if }\langle {\rm wt}(b_1),\alpha^\vee_i\rangle=0, \\ 
\te_i b_1\otimes b_2, & \text{if }\langle {\rm wt}(b_1),\alpha^\vee_i\rangle>0,
\end{cases}
\\
\tf_i(b_1\otimes b_2)=&
\begin{cases}
b_1\otimes \tf_i b_2, & \text{if }\langle {\rm wt}(b_1),\alpha^\vee_i\rangle =0, \\ 
\tf_ib_1\otimes b_2, & \text{if }\langle {\rm wt}(b_1),\alpha^\vee_i\rangle>0,
\end{cases}
\end{split}
\end{equation}}

\item[(3)] 
if $i\in I_{\rm even}$ and $(\e_i,\e_{i+1})=(0,0)$, then
{\allowdisplaybreaks
\begin{equation}\label{eq:tensor product rule for even +}
\begin{split}
&\te_i(b_1\otimes b_2)= \begin{cases}
 b_1 \otimes \te_ib_2, & \text{if $\varphi_i(b_2)\geq\varepsilon_i(b_1)$}, \\ 
\te_ib_1 \otimes  b_2, & \text{if $\varphi_i(b_2)<\varepsilon_i(b_1)$},\\
\end{cases}
\\
&\tf_i(b_1\otimes b_2)=
\begin{cases}
 b_1 \otimes \tf_ib_2, & \text{if $\varphi_i(b_2)>\varepsilon_i(b_1)$}, \\
 \tf_ib_1 \otimes b_2, & \text{if $\varphi_i(b_2)\leq\varepsilon_i(b_1)$}, 
\end{cases}
\end{split}
\end{equation}}
\item[(4)] 
if $i\in I_{\rm even}$ and $(\e_i,\e_{i+1})=(1,1)$, then
{\allowdisplaybreaks
\begin{equation}\label{eq:tensor product rule for even -}
\begin{split}
&\te_i(b_1\otimes b_2)= \begin{cases}
\te_i b_1 \otimes b_2, & \text{if $\varphi_i(b_1)\geq\varepsilon_i(b_2)$}, \\ 
 b_1 \otimes \sigma_i \te_ib_2, & \text{if $\varphi_i(b_1)<\varepsilon_i(b_2)$},\\
\end{cases}
\\
&\tf_i(b_1\otimes b_2)=
\begin{cases}
\tf_ib_1 \otimes  b_2, & \text{if $\varphi_i(b_1)>\varepsilon_i(b_2)$}, \\
b_1 \otimes \sigma_i \tf_i  b_2, & \text{if $\varphi_i(b_1)\leq\varepsilon_i(b_2)$}, 
\end{cases}
\end{split}
\end{equation}}
where $\sigma_i=(-1)^{({\rm wt}(b_1),\alpha_i)}$.
\end{itemize}
\end{prop}
\pf The proof is almost the same as in \cite[Proposition 3.4]{KO}, where the order of tensor product is reversed due to a different convention of comultiplication.
\qed

\begin{rem}\label{rem:natural repn and crystal}
{\rm
Let ${\mc V}=\bigoplus_{i\in \I}\Q(q)v_i$ denote the $\U(\e)$-module, where 
\begin{equation}\label{eq:natural repn}
\omega_i v_j = q_i^{\de_{ij}}v_j,\quad
e_k v_j = \de_{k\, j-1}v_{j-1},\quad f_k v_j =\de_{kj}v_{j+1},
\end{equation}
for $i,j\in \I$ and $k\in I$. It is clear that the pair $\mc{L}=\bigoplus_{i\in \I}A_0 v_i$ and $\mc{B}=\{\,\pm v_i \pmod{q\mc{L}} \,|\,i\in I\,\}$ is a crystal base of $\mc V$.
The crystal structure on $\mc{B}^{\ot \ell}/\{\pm 1\}$ for $\ell\geq 1$ can be described explicitly by Proposition \ref{prop:tensor product rule}, which is the same as in \cite{BKK} or \cite{KO} except that the tensor product order is reversed.  
}
\end{rem}

\section{Schur-Weyl duality and polynomial representations of $\ring\U(\e)$}\label{sec:poly reps}

\subsection{Schur-Weyl duality}
Put $\ring{I}=I\setminus \{0\}$.
Let $\ring{\U(\e)}$ be the $\Q(q)$-subalgebra of $\U(\e)$ generated by $q^h$ and $e_i, f_i$ for $h\in P^\vee$ and $i\in \ring{I}$.

Let us consider ${\mc V}=\bigoplus_{i\in \I}\Q(q)v_i$ in \eqref{eq:natural repn} as a $\ring\U(\e)$-module.
Fix $\ell\ge 2$. Let $\Phi_\ell : \ring\U(\e)\longrightarrow {\rm End}_{\Q(q)}({\mc V}^{\ot \ell})$ denote the action of $\ring\U(\e)$ on ${\mc V}^{\ot \ell}$ via \eqref{eq:comult-1}. Note that $\mc V^{\ot \ell}$ is semisimple (see \cite[Corollary 4.1]{KO}).
  
Assume that $\e_1=0$. We have a $\ring\U(\e)$-linear map $R : {\mc V}^{\ot 2}\longrightarrow {\mc V}^{\ot 2}$ given by
\begin{equation}\label{eq:finite R matrix}
\begin{split}
R(v_i\ot v_j) =
\begin{cases}
q^{-1}q_i^{-1} v_i\ot v_j, & \text{if $i=j$},\\
q^{-1} v_j\ot v_i, & \text{if $i>j$},\\
(q^{-2}-1)v_i\ot v_j + q^{-1} v_j\ot v_i, & \text{if $i<j$},
\end{cases}
\end{split}
\end{equation}
satisfying the Yang-Baxter equation;
\begin{equation*}
R_{12}R_{23}R_{12} = R_{23}R_{12}R_{23},
\end{equation*}
where $R_{ij}$ denotes the map acting as $R$ on the $i$-th and the $j$-th component and the identity elsewhere on ${\mc V}^{\ot 3}$ (cf. \cite{J}).

Let $\mc{H}_\ell(q^{-2})$ be the Iwahori-Hecke algebra of type $A$ over $\Q(q)$ generated by $h_i$ for $i\in \{1,\dots,\ell-1\}$ subject to the relations;
\begin{equation*}
\begin{split}
&(h_i-q^{-2})(h_i+1)=0, \\
& h_ih_j =h_jh_i, \quad\quad\ (|i-j|>1),\\ 
& h_ih_jh_i=h_jh_ih_j \quad (|i-j|=1),
\end{split}
\end{equation*}
for $i, j\in \{1,\dots,\ell-1\}$.
Let $W$ be the symmetric group on $\{1,\dots,\ell\}$ and $s_i=(i\ i+1)$ be the transposition for $1\leq i\leq \ell-1$. For $w\in W$, $\ell(w)$ denote the length of $w$ and let $h(w)$ be the element in $\mc{H}_\ell(q^{-2})$ associated to $w$ such that $h(s_i)=h_i$ for $1\leq i\leq \ell-1$.

We can check that there exists a well-defined action of $\mc{H}_\ell(q^{-2})$ on ${\mc V}^{\ot \ell}$, 
say, $\Psi_\ell : \mc{H}_\ell(q^{-2}) \longrightarrow {\rm End}_{\Q(q)}({\mc V}^{\ot \ell})$,
where $\Psi_\ell(h_i)$ acts as $R$ on the $i$-th and $(i+1)$-th component and the identity elsewhere. Then we have an analogue of Schur-Weyl duality for $\ring\U(\e)$ (cf.~\cite{J}) as follows. The proof is similar to the case when $\e_i=0$ for all $i$. 

\begin{thm}\label{thm:Schur-Weyl duality}
We have
\begin{equation*}
{\rm End}_{\mc{H}_\ell(q^{-2})}({\mc V}^{\ot \ell}) = \Phi_\ell(\ring\U(\e)), \quad 
{\rm End}_{\ring\U(\e)}({\mc V}^{\ot \ell}) = \Psi_\ell({\mc H}_\ell(q^{-2})).
\end{equation*}
\end{thm}
\qed

\subsection{Polynomial representations of $\ring\U(\e)$}

Recall that $M$ is the number of $i$'s with $\e_i=0$ and $N$ is the number of $i$'s with $\e_i=1$ in $\e$.

Let $\cP$ be the set of all partitions. 
A partition $\la=(\la_i)_{i\ge 1}\in \cP$ is called an $(M|N)$-hook partition if $\la_{M+1}\leq N$ (cf.~\cite{BR}). 
We denote the set of all $(M|N)$-hook partitions by $\cP_{M|N}$. For a Young diagram  $\la$, a tableau $T$ obtained by filling $\la$ with letters in $\I$ is called semistandard  if (1) the letters in each row (resp. column) are
weakly increasing from left to right (resp. from top to bottom), (2)
the letters in $\I_0$ (resp. $\I_1$) are strictly increasing in each
column (resp. row).
Let ${SST}_\e(\la)$ be the set of all
semistandard tableaux of shape $\la$. Then $SST_\e(\la)$ is non-empty if and only if $\la\in \cP_{M|N}$. For $T\in SST_\e(\la)$, let $w(T)$ be the word given by reading the entries in $T$ column by column from left to right, and from bottom to top in each column.

For $T\in SST_{\e}(1^r)$ with $r\geq 1$, let $d(T)=\sum_{u<v}d_ud_v$, where $d_u$ is the number of occurrences of $u$ in $T$ for $u\in \I$. 
In general, for a column-semistandard tableau $T$, that is, each column of $T$ is semistandard, we define $d(T) = \sum_{k\geq 1}d(T_k)$, where $T_k$ is the $k$-th column from the left.

We fix $\ell\geq 2$, and let $W$ denote the symmetric group on $\{1,\dots,\ell\}$.
Suppose that $\la\in\cP$ is given with $\sum_{i\geq 1}\la_i=\ell$.
Let $T_+^\la$ be the standard tableau obtained by filling $\la$ with $\{1,\dots,\ell\}$ row by row from top to bottom and from left to right in each row, and let $T_-^\la$ be the tableau obtained by filling $\la$ with $\{1,\dots,\ell\}$ column by column from left to right and from top to bottom in each column.

Let $w_\la\in W$ be such that $w_\la (T_-^\la)=T_+^\la$, where $w_\la (T_-^\la)$ is the tableau obtained by acting $w_\la$ on the letters in $T_-^\la$. Let $W^\la_+$ (resp. $W^\la_-$) be the Young subgroup of $W$ stabilizing the rows (resp. columns) of $T^\la_+$ (resp. $T^\la_-$). 
Then the $q$-deformed Young symmetrizer is given by
\begin{equation}
Y^\la(q) = h(w_\la^{-1})e^\la_+ h(w_\la) e^\la_-,
\end{equation}
(\cite{G}) where
\begin{equation*}
e_+^\la =\sum_{w\in W^\la_+} h(w),\quad 
e_-^\la =\sum_{w\in W^\la_-}(-q^2)^{-\ell(w)} h(w).
\end{equation*}

For $1\leq u<v\leq \ell$, let $W_{uv}=\langle\, s_i \,|\, u\leq i\leq v-1 \rangle$. 
Suppose that $a$ is a letter in $T^\la_-$ such that $a+1$ is located in the same column. We put $C_a=1+h_a$. Then we have 
\begin{align}
Y^\la(q) C_a &=0,\label{eq:column relation}
\end{align}
Next, suppose that $a$ is a letter in $T^\la_-$, where there is another letter $d$ to the right. Let $b$ be the letter at the bottom of column where $a$ is placed, and $c=b+1$ the letter at the top of the column where $d$ is placed. Let $\mc{G}^\la_a$ be the set of minimal length right coset representatives of $W_{ab}\times W_{cd}$ in $W_{ad}$. We define the {\em Garnir element} at $a$ to be
\begin{equation}
G^\la_a = \sum_{w\in \mc{G}^\la_a} (-q^2)^{\ell(w)} h(w).
\end{equation} 
The collection of boxes in the Young diagram $\la$ corresponding to the letters from $a$ to $d$ in $T^\la_-$ is called a {\em Garnir belt} at $a$.
Then we have the following relations \cite[(15)]{BKW};
\begin{align}
Y^\la(q) G^\la_a &=0 \label{eq:Garnir relation}.
\end{align}

Let $T$ be a tableau  of shape $\la$ with letters in $\I$, and let $T(i)$ be the letter in $T$ at the position corresponding to $i$ in $T^\la_-$ for $1\leq i\leq \ell$.
Let 
\begin{equation*}
v_T = Y^\la(q)\left(v_{T(1)}\ot \dots \ot v_{T(\ell)}\right).
\end{equation*}
For $\sigma\in W$, let $T^\sigma$ be the tableau given by replacing $T(i)$ with $T(\sigma(i))$ for $1\leq i\leq \ell$.

Let $a$ be a letter in $T^\la_-$ with $d$ to the right in the same row and with $b,c$ as above.
Let $w_0$ be the longest element in $W_{ab}\times W_{cd}$, 
and let $\tilde{\mc{G}}^\la_a = w_0{\mc{G}}^\la_a w_0$.
Let $u_1,\dots, u_s$ and $u_{s+1},\dots, u_{r+s}$ be the letters in $T$ corresponding to $c,\dots,d$ and $a,\dots,b$ in $T^\la_-$, respectively.
Then we may identify $\sigma\in \tilde{\mc{G}}^\la_a$ with a permutation on $\{1,\dots,r+s\}$ 
satisfying $\sigma(1)<\dots<\sigma(s)$ and $\sigma(s+1)<\dots<\sigma(s+r)$ so that $T^{\sigma}$ is the tableau obtained from $T$ by replacing $u_i$'s with $u_{\sigma(i)}$'s for $1\leq i\leq r+s$. 
With this identification, we let $\tilde{\ell}(\sigma)$ be the length of $\sigma$ as a permutation on $\{1,\dots,r+s\}$, and put
\begin{equation*}\label{eq:XY for sigma}
\begin{split}
X_{\sigma}&=\{\,i\,|\,1\leq i\leq s,\ s+1\leq \sigma^{-1}(i)\leq s+r\,\},\\
Y_{\sigma}&=\{\,j\,|\,s+1\leq j\leq s+r,\ 1\leq \sigma^{-1}(j)\leq r\,\}.
\end{split}
\end{equation*}


\begin{lem}\label{lem:Garnir relation}
Suppose that $T$ is column-semistandard such that either $T(a)=T(d)\in \I_1$ or  $T(a)>T(d)$. Then under the above hypothesis, we have
\begin{equation*}
v_T= - \sum_{\sigma\in \tilde{\mc{G}}^\la_a,\sigma \ne 1}
(-q)^{{\tilde{\ell}(\sigma)}+m(\sigma,T)}v_{T^{\sigma}},
\end{equation*}
where
\begin{equation*}
\begin{split}
m(\sigma,T) =
& - \big\vert\{\,(i,j)\,|\,1\leq i<j \leq s,\ {i \notin X_\sigma, j\in X_\sigma},\ u_i=u_j\,\} \big\vert \\
& - \big\vert\{\,(k,l)\,|\ s+1 \leq k<l \leq s+r,\ {k\in Y_\sigma, l \notin Y_\sigma},\ u_k=u_l\,\} \big\vert \\
& + \big\vert\{\,(x,y)\,|\,1 \leq x\leq s,\ s+1\leq y\leq s+r,\ {\text{$x\in X_\sigma$ or $y\in Y_\sigma$}},\ u_x=u_y\,\} \big\vert. 
\end{split}
\end{equation*}
\end{lem}

\pf We have $v_T = Y^\la(q) v$, where $v=\left(v_{T(1)}\ot \dots \ot v_{T(\ell)}\right)$. 
Following the above notations, we have $v=v'\ot v_{u_{1+s}} \ot \cdots \ot v_{u_{r+s}} \otimes v_{u_{1}} \otimes \cdots v_{u_{s}} \ot v''$. 
Note that 
\begin{equation*}
v_{T^{w_0}} = Y^\la(q)\left( 
v'\ot v_{u_{r+s}} \ot \cdots \ot v_{u_{1+s}} \otimes v_{u_{s}} \otimes \cdots v_{u_{1}} \ot v''
\right),
\end{equation*}
where $u_{r+s}\ge \dots u_{1+s}=T(a)\ge u_s=T(d) \ge \dots \ge u_1$.

For $w\in {\mc{G}}^\la_a$, we have by \eqref{eq:finite R matrix} and \eqref{eq:column relation}
\begin{equation}\label{columngarnir}
\begin{split}
& h(w)\left( 
v'\ot v_{u_{r+s}} \ot \cdots \ot v_{u_{1+s}} \otimes v_{u_{s}} \otimes \cdots v_{u_{1}} \ot v''
\right)\\
=&q^{-\ell(w)}(-q)^{m(\sigma,T^{w_0})} 
\left( 
v'\ot v_{u_{\sigma(r+s)}} \ot \cdots \ot v_{u_{\sigma(1+s)}} \otimes v_{u_{\sigma(s)}} \otimes \cdots \otimes v_{u_{\sigma(1)}} \ot v''
\right),
\end{split}
\end{equation}
where $\sigma$ is the permutation on $\{1,\dots,r+s\}$ corresponding to $w_0ww_0$ and \begin{equation*}
m(\sigma,T^{w_0}) = \big\vert \{\,(i,j)\,|\,i<j,\ \sigma^{-1}(i)< \sigma^{-1}(j),\ u_i=u_j \,\} \big\vert.
\end{equation*}
Hence it follows from \eqref{eq:Garnir relation} and that \eqref{columngarnir}
\begin{equation*}
\begin{split}
0&=Y^{\la}(q)G^{\la}_a(q)
\left( 
v'\ot v_{u_{r+s}} \ot \cdots \ot v_{u_{1+s}} \otimes v_{u_{s}} \otimes \cdots v_{u_{1}} \ot v''
\right) \\
&=Y^{\la}(q)\sum_{w\in \mc{G}^\la_a}(-q^2)^{\ell(w)}h(w)
\left( 
v'\ot v_{u_{r+s}} \ot \cdots \ot v_{u_{1+s}} \otimes v_{u_{s}} \otimes \cdots v_{u_{1}} \ot v''
\right)\\
&=Y^{\la}(q)\sum_{w\in \mc{G}^\la_a}(-q)^{\ell(w)+m(\sigma,T^{w_0})}
\left( 
v'\ot v_{u_{\sigma(r+s)}} \ot \cdots \ot v_{u_{\sigma(1+s)}} \otimes v_{u_{\sigma(s)}} \otimes \cdots v_{u_{\sigma(1)}} \ot v''
\right)\\
&=\sum_{w\in \mc{G}^\la_a}(-q)^{\ell(w)+m(\sigma,T^{w_0})}v_{T^{w_0w}}
=\sum_{\sigma\in \tilde{\mc{G}}^\la_a}(-q)^{\tilde\ell(\sigma)+ m(\sigma,T^{w_0})}v_{T^{\sigma w_0}}.
\end{split}
\end{equation*}

We have
\begin{equation}\label{eq:aux0}
\sum_{\sigma\in \tilde{\mc{G}}^\la_a}
(-q)^{\tilde\ell({\sigma})+m(\sigma, T^{w_0})}v_{(T^{\sigma})^{w_0}}=0.
\end{equation}

For $\sigma\in \tilde{\mc{G}}^\la_a$,
let $U^{\sigma}$ be the subtableau of $T^{\sigma}$ corresponding to the Garnir belt at $a$, where $U=U^{\rm id}$. 
We define $d_a(T^\sigma)$ in the same way as in $d(T)$ only by using the letters in $U^\sigma$.
Let {$l_p>\dots>l_1\geq r_q>\dots>r_1$} be the distinct letters appearing in $U$, where
$l_i$ and $r_j$ are located in the left and right columns of $U$, respectively.

Let $m_i$ (resp. $n_j$) be the number of occurrences of $l_i$'s (resp. $r_j$'s) in $U$, which remain in the same column after applying $\sigma$. Let $m'_i$ (resp. $n'_i$) be the number of $l_i$'s (resp. $r_j$'s) which are placed on the right (resp. left) column of $U^\sigma$ after applying $\sigma$ to $U$. 
Note that $\sum_{i}m'_i=\sum_{j}n'_j$

{\em Case 1}. Suppose that $l_1\neq r_q$.
We have
\begin{equation}\label{eq:aux1}
\begin{split}
d_a(T) &= \sum_{1\le i<j\le p}(m_i+m'_i)(m_j+m'_j) + \sum_{1\le k<l\le q}(n_k+n'_k)(n_l+n'_l),
\end{split}
\end{equation}
while 
\begin{equation}\label{eq:aux2}
\begin{split}
&d_a(T^\sigma) \\
&=
\sum_{1\le i<j\le p}(m_im_j + m'_im'_j) + \sum_{i,k}m_in'_k 
+\sum_{1\le k<l\le q}(n_kn_l + n'_kn'_l) + \sum_{j,l}m'_jn_l\\
&=\sum_{1\le i<j\le p}(m_im_j + m'_im'_j) + \sum_{i}m_i \sum_k n'_k
+\sum_{1\le k<l\le q}(n_kn_l + n'_kn'_l) + \sum_{j}m'_j\sum_{l}n_l.\\
\end{split}
\end{equation}

Since we have 
\begin{equation*}
\begin{split}
m(\sigma,T^{w_0}) = 
& \quad\ \big\vert\{\,(i,j)\,|\ 1\leq i<j\leq s,\ {i\in X_{\sigma}, j\not\in X_{\sigma}},\ u_i=u_j\,\} \big\vert \\
& + \big\vert\{\,(k,l)\,|\,s+1\leq k<l \leq s+r,\ {k\not\in Y_{\sigma}, l\in Y_{\sigma}},\ u_k=u_l\,\} \big\vert, \\
m(\sigma,T) =
& - \big\vert\{\,(i,j)\,|\,1\leq i<j \leq s,\ {i \notin X_\sigma, j\in X_\sigma},\ u_i=u_j\,\} \big\vert \\
& - \big\vert\{\,(k,l)\,|\ s+1 \leq k<l \leq s+r,\ {k\in Y_\sigma, l \notin Y_\sigma},\ u_k=u_l\,\} \big\vert,
\end{split}
\end{equation*}
one can check easily that
{
\begin{equation}\label{eq:aux3}
m(\sigma,T^{w_0}) - m(\sigma,T) = \sum_{1\le i\le p}m_im'_i +\sum_{1\le j\le q} n_jn'_j.
\end{equation}}
By \eqref{eq:aux1}, \eqref{eq:aux2}, and \eqref{eq:aux3}, we have
\begin{equation}\label{eq:aux4}
d_a(T) - d_a(T^\sigma) =  m(\sigma,T) - m(\sigma,T^{w_0}).
\end{equation}
By \eqref{eq:column relation}, \eqref{eq:aux0} and \eqref{eq:aux4}, we have
\begin{equation*}\label{eq:aux5}
\begin{split}
0=&
\sum_{\sigma\in \tilde{\mc{G}}^\la_a}
(-q)^{\tilde\ell({\sigma})+m(\sigma, T^{w_0})}v_{(T^{\sigma})^{w_0}}
=
\sum_{\sigma\in \tilde{\mc{G}}^\la_a}
(-q)^{\tilde\ell({\sigma})+m(\sigma, T^{w_0}) - d_a(T^{\sigma})}v_{T^{\sigma}}\\
=&
\sum_{\sigma\in \tilde{\mc{G}}^\la_a}
(-q)^{\tilde\ell({\sigma})+m(\sigma, T) - d_a(T)}v_{T^{\sigma}}
=(-q)^{- d_a(T)}\sum_{\sigma\in \tilde{\mc{G}}^\la_a}
(-q)^{\tilde\ell({\sigma})+m(\sigma, T)}v_{T^{\sigma}}.
\end{split}
\end{equation*}
This proves the identity in the lemma.

{\em Case 2}. Suppose that Suppose that $l_1 = r_q$.
In this case, $d_a(T)$ is the same as in {\em Case 1}, and 
\begin{equation*}
d_a(T)-d_a(T^\sigma) = - \sum_{1\le i\le p}m_im'_i - \sum_{1\le j\le q} n_jn'_j + m_pn'_1 + m'_p n_1.
\end{equation*}
Note that
\begin{equation*}
\begin{split}
m(\sigma,T^{w_0}) = 
& \quad\ \big\vert\{\,(i,j)\,|\ 1\leq i<j\leq s,\ {i\in X_{\sigma}, j\not\in X_{\sigma}},\ u_i=u_j\,\} \big\vert  \\
& + \big\vert\{\,(k,l)\,|\,s+1\leq k<l \leq s+r,\ {k\not\in Y_{\sigma}, l\in Y_{\sigma}},\ u_k=u_l\,\} \big\vert  \\
&+ \big\vert\{\,(x,y)\,|\,1 \leq x\leq s,\ s+1\leq y\leq s+r,\ x\in X_\sigma, y\in Y_\sigma, u_x=u_y\,\} \big\vert ,
\end{split}
\end{equation*}
where the last summand is equal to $m'_pn'_1$.
By similar arguments as in \eqref{eq:aux3}, we have
\begin{equation*}\label{eq:aux6}
d_a(T) - d_a(T^\sigma) =  m(\sigma,T) - m(\sigma,T^{w_0}).
\end{equation*}
This also proves the identity in the lemma as in \eqref{eq:aux4}.
\qed
\vskip 2mm

For $\la\in \cP_{M|N}$ with $\sum_{i}\la_i=\ell$, let 
\begin{equation}\label{eq:Ve(la)}
V_\e(\la)=\sum_{T\in SST_\e(\la)}\Q(q)v_T.
\end{equation} 
Let $H_\la$ be the tableau in $SST_\e(\la)$, which is defined inductively as follows:
\begin{itemize}
\item[(1)] Fill the first row (resp. column) of $\la$ with $1$ if $\e_1=0$ (resp. $\e_1=1$).

\item[(2)] Suppose that we have filled a subdiagram of $\la$ from $1$ to $i$. 
Then fill the first row (resp. column) of the remaining diagram with $i+1$ if $\e_{i+1}=0$ (resp. $\e_{i+1}=1$).
\end{itemize}

\begin{ex}{\rm
Suppose that $n=5$, $\e=(0,1,1,0,0)$ and $\lambda=(6,5,4,2,1)$. In this case, we have
\begin{equation*}
H_{\lambda}= 
\raisebox{3ex}{
\ytableausetup {mathmode, boxsize=1.0em} 
\begin{ytableau}
 {1} & {1}  & {1}& {1}& {1}& {1}\\
  {2} &  {3}& {4}& {4}&{4}  \\
  {2} &  {3}& {5}& {5}    \\
 {2} & {3} \\
{2}  \\
\end{ytableau}}
\end{equation*}

}
\end{ex}

\begin{prop}\label{prop:poly repn}
We have the following.
\begin{itemize}
\item[(1)] $V_\e(\la)$ is a $\ring\U(\e)$-submodule of $\V^{\ot\ell}$.

\item[(2)] $V_\e(\la)$ is an irreducible $\ring\U(\e)$-module with basis $\{\,v_T\,|\,T\in SST_\e(\la)\,\}$.

\item[(3)] $v_{H_\la}$ is a highest weight vector in $V_\e(\la)$.
\end{itemize}
\end{prop}
\pf (1) It is clear that $V_\e(\la)$ is invariant under $q^{h}$ for $h\in P^\vee$. 
It suffices to check $f_i V_\e(\la)\subset V_\e(\la)$ for $i\in \ring{I}$ since the proof for $e_i$ is the same.
The proof is similar to the case when $\e=(0,\dots,0)$ (cf.~\cite{Ful}). For column-semistandard tableaux $U$ and $V$ of shape $\la$, we define $U<V$ if there exists $1\le k\le \ell$ such that $U(k)<V(k)$ and $U(k')=V(k')$ for $k<k'\le \ell$. 

Suppose that $T\in SST_\e(\la)$ is given. By \eqref{eq:comult-1}, $f_i v_T$ is a $\Q(q)$-linear combination of $v_{T'}$'s, where we may assume that $T$ is column-semistandard by \eqref{eq:column relation}. If such $T'$ is not semistandard, then we may apply Lemma \ref{lem:Garnir relation} to $T'$ so that $v_{T'}$ is a linear combination of $T''$'s which is column-semistandard and $T'<T''$. Repeating this process finitely many times, we conclude that $f_iT$ is a linear combination of $v_S$'s for some $S\in SST_\e(\la)$. Therefore, have $f_iV_\e(\la)\subset V_\e(\la)$.

(2) Since $V_\e(\la) = Y^\la(q)\V^{\ot \ell}$ and $Y^\la(q)$ is a primitive idempotent up to scalar multiplication \cite{G}, it follows from Theorem \ref{thm:Schur-Weyl duality} that $V_\e(\la)$ is an irreducible $\ring\U(\e)$-module. Recall that the dimension of the irreducible $\mc{H}_\ell(q^{-2})$-module $S^\la$ generated by $Y^\la(q)$ is the number of standard tableaux of shape $\la$.
We may have an analogue of the Robinson-Schensted type correspondence, which is a bijection from the set of words of length $\ell$ with letters in $\I$ to the set of pair of standard tableau and semistandard tableau of shape $\la$ (cf.~\cite{BR}). Comparing the dimensions of $\V^{\ot \ell}$ and its decomposition into $\mc{H}(q^{-2})\ot \ring\U(\e)$-module $S^\la\ot V_\e(\la)$, we conclude that $\dim_{\Q(q)}V_\e(\la)$ is equal to $|SST_\e(\la)|$, and hence $\{\,v_T\,|\,T\in SST_\e(\la)\,\}$ is a linear basis of $V_\e(\la)$. 

(3) The character of $V_\e(\la)$ is equal to that of polynomial representations of the general linear Lie superalgebra $\mf{gl}_{M|N}$ corresponding to $\la\in \cP_{M|N}$, and ${\rm wt}(v_{H_\la})$ is maximal \cite[Theorem 2.55]{CW}. This implies that $e_i v_{H_\la}=0$ for all $i\in \ring{I}$ and hence $v_{H_\la}$ is a highest weight vector.
\qed

\begin{rem}\label{rem:character of poly}
{\rm
The character of $V_\e(\la)$ for $\la\in \cP_{M|N}$ is called a hook Schur polynomial \cite{BR}, which depends only on $\e$ up to permutations. 
The tensor product of two polynomial representations is completely reducible and the multiplicity of each irreducible component is given by usual Littlewood-Richardson coefficient. 
}
\end{rem}

\subsection{Crystal base of $V_\e(\la)$}
Let $\la\in \cP_{M|N}$ be given. 
We may define an $\ring I$-colored oriented graph structure by identifying $T$ with $w(T)^{\rm rev}$, the reverse word of $w(T)$.

Let
\begin{equation}\label{eq:LB for lambda}
\begin{split}
L_\e(\la) &=\bigoplus_{T\in SST_\e(\la)}A_0 v_T^*,\\ 
B_\e(\la) &=\{\,\pm v_T^*\!\!\!\!\pmod{qL_\e(\la)}\,|\,T\in SST_\e(\la)\,\},
\end{split}
\end{equation}
where $v_T^*=q^{-d(T)}v_T$ for $T\in SST_\e(\la)$.

\begin{lem}\label{lem:crystal base of column and row}
When $\la=(1^r)$ or $(r)$ for $r\geq 1$,  $(L_\e(\la),B_\e(\la))$ is a crystal base of $V_\e(\la)$, and the crystal $B_\e(\la)/\{\pm 1\}$ is isomorphic to $SST_\e(\la)$.
\end{lem}
\pf The proof is similar to that of \cite[Proposition 3.3]{KO}.
\qed

\begin{prop}
Suppose that $\e=\e_{M|N}$. 
For $\la\in \cP_{M|N}$, 
$(L_\e(\la),B_\e(\la))$ is a crystal base of $V_\e(\la)$.
\end{prop}
\pf The proof is similar to that of \cite[Theorem 4.4]{LT}.
Let $(\mc{L}(\la),\mc{B}(\la))$ be given by
\begin{equation*}\label{eq:(L,B)}
\begin{split}
\mc{L}(\la)&=\sum_{r\geq 0,\, i_1,\ldots,i_r\in \ring I}A_0 \td{x}_{i_1}\cdots\td{x}_{i_r}v_\la, \\
\mc{B}(\la)&=\{\,\pm \,{\td{x}_{i_1}\cdots\td{x}_{i_r}v_\la}\!\!\! \mod{q L(\la)}\,|\,r\geq 0, i_1,\ldots,i_r\in \ring I\,\}\setminus\{0\},
\end{split}
\end{equation*}
where $v_\la$ is a highest weight vector in $V_\e(\la)$ and $x=e, f$ for each $i_k$.  
Following the same arguments in \cite{BKK}, it is shown in \cite{KO} that $(\mc{L}(\la),\mc{B}(\la))$ is a crystal base of $V_\e(\la)$. The crystal $\mc{B}(\la)/\{\pm 1\}$ is equal to $SST_\e(\la)$ which is connected.

Let $\mu=(\mu_1,\dots,\mu_r)=\la'$ be the conjugate partition of $\la$, and 
\begin{equation*}
V_\e^\mu = V_\e((1^{\mu_1}))\ot \dots \ot V_\e((1^{\mu_r})).
\end{equation*}
Let $I^\mu_\e$ be the subspace of $V_\e^\mu$ spanned by the vectors induced from the relation \eqref{eq:Garnir relation}, which includes the relations in Lemma \ref{lem:Garnir relation}.
Since $I^\mu_\e$ is a $\ring\U(\e)$-submodule, the quotient $V_\e^\mu/I^\mu_\e$ is isomorphic to $V_\e(\la)$ by Proposition \ref{prop:poly repn}. 
So we have a well-defined $\ring\U(\e)$-linear map
\begin{equation*}
\pi^\mu : V^\mu_\e \longrightarrow V_\e(\la)
\end{equation*}
given by $\pi^\mu(v_{T_1}\ot \dots \ot v_{T_r})=v_T$ where $T$ is the column semistandard tableau whose $i$-th column (from the left) is $T_i$ for $1\leq i\leq r$.
Since the decomposition of $V^{\mu}$ is equal to the usual Pieri rule of Schur functions, it has exactly one component isomorphic to $V_\e(\la)$. Hence $\pi^\mu$ is equal to the projection onto $V_\e(\la)$ up to scalar multiplication. 

Let $L_\e^\mu = L_\e((1^{\mu_1}))\ot \dots \ot L_\e((1^{\mu_r}))$ be the crystal lattice of $V^\mu_\e$. By \cite[Theorem 4.14]{KO}, $\pi^\mu(L_\e^\mu)$ is a crystal lattice of $V_\e(\la)$ whose ${\rm wt}(H_{\la})$-weight space is equal to $A_0v^*_{H_\la}$.
Since the crystal of $V_\e(\la)$ is connected, we conclude that $\{\,v^*_T\,|\,T\in SST_\e(\la)\,\}$ is an $A_0$-basis of $\pi^\mu(L_\e^\mu)$ which is equal to $L_\e(\la)$.
\qed

\begin{rem}\label{rem:crystal base of poly}
{\rm
For arbitrary $\e$, the $\ring I$-colored oriented graph $SST_\e(\la)$ is not in general connected (see \cite{K09} for more details). Furthermore, it is not known yet whether $V_\e(\la)$ has a crystal base for any $\la\in \cP_{M|N}$. We expect that $(L_\e(\la),B_\e(\la))$ in \eqref{eq:LB for lambda} is a crystal base of $V_\e(\la)$.
}
\end{rem}

\section{$R$ matrix for finite-dimensional $\U(\e)$-modules}\label{sec:R matrix}

\subsection{Finite-dimensional $\U(\e)$-modules of fundamental type}

Let $\Z_+$ be the set of non-negative integers.
Let
\begin{equation*}
\Z^n_+(\e)=\{\,{\bf m}=(m_1,\ldots,m_n)\,|\,\text{$m_i\in \Z_+$ if $\e_i=0$, $m_i\in \{0,1\}$ if $\e_i=1$, ($i\in \I$)}\,\}.
\end{equation*}
For ${\bf m}\in \Z_+^n(\e)$, let $|{\bf m}|=m_1+\dots+m_n$.
For $i\in \I$, put $\be_i=(0,\cdots, 1,\cdots, 0)$ where $1$ appears only in the $i$-th component. 

For $s\in \Z_{+}$, let 
\begin{equation*}
\W_{s,\e} = \bigoplus_{\substack{{\bf m}\in \Z^n_+(\e), |{\bf m}|=s}}\Q(q)|{\bf m}\rangle
\end{equation*}
be the $\Q(q)$-vector space spanned by $|{\bf m}\rangle$ for ${\bf m}\in \Z^n_+(\e)$ with $|{\bf m}|=s$.

For a parameter $x\in \Q(q)$, we denote by $\W_{s,\e}(x)$ a $\U(\e)$-module $V$, where $V=\W_{s,\e}$ as a $\Q(q)$-space and the actions of $e_i, f_i, \omega_j$ are given by  
\begin{equation*}\label{eq:W_l}
\begin{split}
e_i |{\bf m}\rangle &= 
\begin{cases}
x^{\delta_{i,0}}{q^{m_{i+1}-m_{i}-1}}[m_{i+1}]|{\bf m} + \be_{i} -\be_{i+1} \rangle, & \text{if ${\bf m} + \be_{i} -\be_{i+1}\in\Z_+^n(\e)$},\\
0,& \text{otherwise},
\end{cases}
\\
f_i |{\bf m}\rangle &= 
\begin{cases}
x^{-\delta_{i,0}}{q^{m_{i}-m_{i+1}-1}}[m_{i}]|{\bf m} - \be_{i} + \be_{i+1} \rangle, & \text{if ${\bf m} - \be_{i} + \be_{i+1}\in\Z_+^n(\e)$},\\
0,& \text{otherwise},
\end{cases}\\
\omega_j |{\bf m}\rangle &= q_j^{m_j} |{\bf m} \rangle,
\end{split}
\end{equation*}
for $i\in I$, $j\in \I$, and ${\bf m}=(m_1,\ldots,m_n)\in \Z^n_+(\e)$. 
Here we understand $\be_0=\be_n$. 

\begin{rem}{\rm
We may identify $\W_{s,\e}(x)$ with $V_\e((s))$ \eqref{eq:Ve(la)} as a $\ring\U(\e)$-module, where $|{\bf m}\rangle$ corresponds to $v_T$, where $T$ is the tableau of shape $(s)$ with $m_i$ the number of occurrences of $i$ in $T$ ($i\in \I$). Also the map 
\begin{equation}\label{eq:iso phi}
\phi(|{\bf m}\rangle)= q^{-\sum_{i<j}m_im_j}|{\bf m}\rangle
\end{equation}
gives an isomorphism of $\U(\e)$-modules from $\W_{s,\e}(x)$ to itself with another $\U(\e)$-action defined in \cite[(2.15)]{KO}.}
\end{rem}

Let us  regard $\W_{s,\e}=\W_{s,\e}(1)$ and set
\begin{equation}\label{eq:crystal base of W(r)}
\begin{split}
\mc{L}_{s,\e} &=\bigoplus_{{\bf m}\in \Z_+^n(\e), |{\bf m}|=s}A_0 | {\bf m} \rangle,\quad
\mc{B}_{s,\e} =\{\,\pm |{\bf m}\rangle \!\!\! \pmod{q\mc{L}_{s,\e}}\,|\ {\bf m}\in \Z_+^n(\e), |{\bf m}|=s\,\}.
\end{split}
\end{equation}

\begin{prop}  
For $s\in \Z_{+}$, the pair $(\mc{L}_{s,\e},\mc{B}_{s,\e})$ is a crystal base of $\W_{s,\e}$, where the crystal $\mc{B}_{s,\e}/\{\pm 1\}$ is connected.
\end{prop}
\pf It follows from the same arguments as in Lemma \ref{lem:crystal base of column and row} that $(\mc{L}_{s,\e},\mc{B}_{s,\e})$ is a crystal base of $\W_{s,\e}$. The crystal $SST_\e((s))$ is connected with highest element $H_{(s)}$. Since the crystal $\mc{B}_{s,\e}/\{\pm 1\}$ of $\W_{s,\e}$ is equal to $SST_\e((s))$ as an $\ring I$-colored graph, $\mc{B}_{s,\e}/\{\pm 1\}$ is connected as an $I$-colored oriented graph.
\qed

\subsection{Subalgebra $\U(\e')$}\label{subsec:subalg e'}
Suppose that {$n\geq 4$} and let $\e=(\e_1,\dots,\e_n)$ be given.
Let $\e'=(\e'_1,\dots,\e'_{n-1})$ be the sequence obtained from $\e$ by removing $\e_i$ for some $i\in \I$. {We further assume that $\e'$ is homogeneous when $n=4$, that is, $\e'=(000)$ or $(111)$.}

Put $I'=\{0,1,\cdots,n-2\}.$ 
Let us denote by $\omega'_l$, $e'_j$, and $f'_j$ the generators of $\U(\e')$ for $1 \leq l \leq n-1 $ and $j \in I'$, where $k'_j=\omega'_j(\omega'_{j+1})^{-1}$.
Let us  define $K_j, E_j, F_j$ for $j \in I'$ as follows:

{\em Case 1}. Assume that $2\leq i\leq n-1$. 
For $j\in I'$, put 
\begin{equation}\label{eq:generators for e'}
\begin{split}
K_j &=
\begin{cases}
k_j, & \text{if $j\leq i-2$},\\
k_{i-1}k_{i}, & \text{if $j= i-1$},\\
k_{j+1}, & \text{if $j\geq i$},
\end{cases}\\
E_j &=
\begin{cases}
e_j, & \text{if $j\leq i-2$},\\
[e_{i-1},e_{i}]_{D_{i-1\,i}}, & \text{if $j= i-1$},\\
e_{j+1}, & \text{if $j\geq i$},
\end{cases}\quad
F_j =
\begin{cases}
f_j, & \text{if $j\leq i-2$},\\
[f_{i},f_{i-1}]_{D_{i-1\,i}^{-1}}, & \text{if $j= i-1$},\\
f_{j+1}, & \text{if $j\geq i$}.
\end{cases}
\end{split}
\end{equation}

{\em Case 2}. Assume that $i=n$. 
For $j\in I'$, put 
\begin{equation}\label{eq:generators for e'-2}
\begin{split}
K_j &=
\begin{cases}
k_j, & \text{if $j\neq 0$},\\
k_{n-1}k_{0}, & \text{if $j=0$},\\
\end{cases}\\
E_j &=
\begin{cases}
e_j, & \text{if $j\neq 0$},\\
[e_{n-1},e_{0}]_{D_{n-1\,0}}, & \text{if $j=0$},\\
\end{cases}\quad
F_j =
\begin{cases}
f_j, & \text{if $j\neq 0$},\\
[f_{0},f_{n-1}]_{D_{n-1\,0}^{-1}}, & \text{if $j=0$}.\\
\end{cases}
\end{split}
\end{equation}

{\em Case 3}. Assume that $i=1$. 
For $j\in I'$, put 
\begin{equation}\label{eq:generators for e'-3}
\begin{split}
K_j &=
\begin{cases}
k_{0}k_{1}, & \text{if $j=0$},\\
k_{j+1}, & \text{if $j\neq 0$},
\end{cases}\\
E_j &=
\begin{cases}
[e_{0},e_{1}]_{D_{0\,1}}, & \text{if $j=0$},\\
e_{j+1}, & \text{if $j\neq 0$},
\end{cases}\quad
F_j =
\begin{cases}
[f_{1},f_{0}]_{D_{0\,1}^{-1}}, & \text{if $j=0$},\\
f_{j+1}, & \text{if $j\neq 0$}.
\end{cases}
\end{split}
\end{equation}

\begin{thm}\label{thm:folding homomorphism}
There exists a homomorphism of $\Q(q)$-algebras $\phi : \U(\e') \longrightarrow \U(\e)$ such that
\begin{equation*}
\phi(k'_j)=K_j,\quad \phi(e'_j)=E_j,\quad \phi(f'_j)=F_j\quad (j\in I').
\end{equation*}
\end{thm}
\pf 
Let us prove {\em Case 1} since the the proof of the other cases are similar.
Let $\td\e=(\td\e_1,\dots\td\e_n)$ be the sequence obtained from $\e$ by exchanging $\e_i$ and $\e_{i+1}$, and let $\tau_i : \U(\e) \longrightarrow \U(\td\e)$ be the isomorphism in Theorem \ref{thm:reflection}.

Put $\Omega_j=\omega_j$ for $1\leq j \leq i-1$ and $\omega_{j+1}$ for $i\leq j \leq n-1$, and let $\phi(\omega'_j)=\Omega_j, \phi(e'_j)=E_j,$ and $\phi(f'_j)=F_j$ for $j=1,\cdots, n-1$. Let us check that $\Omega_j, E_j, F_j$ satisfy the relations in Definition \ref{def:U(e)}. 
Note that $D_{i-1 i}=q_i^{-1}$.

First, the relations \eqref{eq:Weyl-rel-1} and \eqref{eq:Weyl-rel-2} are trivial. Let us check that \eqref{eq:Weyl-rel-3} holds. Let $E_j$ and $F_l$ be given for $j,l \in I'$. If $j\neq l$ or $j=l\neq i-1$, then it is clear. When $j=l=i-1$, we have $\tau_i^{-1}(e_{i-1})=[e_{i-1},e_i]_{D_{i-1\,i}}=E_{i-1}$, $\tau^{-1}(f_{i-1})=F_{i-1}$, and $\tau^{-1}(k_{i-1})=K_{i-1}$. Hence \eqref{eq:Weyl-rel-3} holds. We can check the relation \eqref{eq:Weyl-rel-4} by the same argument.

Next, consider the relations \eqref{eq:Serre-rel-1}. The first one is immediate. 
So it is enough to show the second one. We may only consider four non-trivial cases when the pair of relevant indices in $I'$ are $(i-2,i-1), (i-1,i-2),(i-1,i),(i,i-1)$ with the first index in the pair in $I'_{\rm even}$. In case of $(i-2,i-1)$, we have 
\begin{equation*}
\begin{split}
&E_{i-2}^2E_{i-1}-(-1)^{\epsilon_{i-2}}[2]E_{i-2}E_{i-1}E_{i-2}+E_{i-1}E_{i-2}^2\\
&{=e_{i-2}^2e_{i-1}e_{i}-q_i^{-1}e_{i-2}^2e_{i}e_{i-1}-(q_{i-1}+q_{i-1}^{-1})e_{i-2}e_{i-1}e_ie_{i-2}}\\
&\quad {+(q_{i-1}+q_{i-1}^{-1})q_i^{-1}e_{i-2}e_ie_{i-1}e_{i-2}+e_{i-1}e_ie_{i-2}^2-q_i^{-1}e_ie_{i-1}e_{i-2}^2.}
\end{split}
\end{equation*}
which is zero, since $e_{i-2}^2e_{i-1}+e_{i-1}e_{i-2}^2=(q_{i-1}+q_{i-1}^{-1})e_{i-2}e_{i-1}e_{i-2}$ and hence
{\begin{equation*}
\begin{split}
&e_{i-2}^2e_{i-1}e_i-(q_{i-1}+q_{i-1}^{-1})e_{i-2}e_{i-1}e_ie_{i-2}+e_{i-1}e_ie_{i-2}^2=0,\\
&-q_i^{-1}e_{i-2}^2e_ie_{i-1}+q_i^{-1}(q_{i-1}+q_{i-1}^{-1})e_{i-2}e_ie_{i-1}e_{i-2}-q_i^{-1}e_ie_{i-1}e_{i-2}^2=0. 
\end{split}
\end{equation*}}
The proof for $(i,i-1)$ is the same. In case of $(i-1,i-2)$ and $(i-1,i)$, the proof reduces to the case of $(i-2,i-1)$ or $(i,i-1)$ by applying $\tau_i$ to $E_l$'s for $l= i-2,i-1,i$.

Finally let us check the relation \eqref{eq:Serre-rel-2}. We may only consider the cases when the relevant triple of indices in $I'$ are $(i-3,i-2,i-1), (i-2,i-1,i), (i-1,i,i+1)$ with the index in the middle in $I'_{\rm odd}$. In case of $(i-1,i,i+1)$ and $i \in I'_{\rm odd}$, we have
\begin{equation*}
\begin{split}
&E_{i}E_{i-1}E_{i}E_{i+1}  
- E_{i}E_{i+1}E_{i}E_{i-1} + E_{i+1}E_{i}E_{i-1}E_{i} \\  
&\hskip 2.2cm - E_{i-1}E_{i}E_{i+1}E_{i} 
+ (-1)^{\e_i}[2]E_{i}E_{i-1}E_{i+1}E_{i} \\
&= e_{i+1}(e_{i-1}e_i-q_{i}^{-1}e_ie_{i-1})e_{i+1}e_{i+2}  
- e_{i+1}e_{i+2}e_{i+1}(e_{i-1}e_i-q_{i}^{-1}e_ie_{i-1}) \\
&\quad + e_{i+2}e_{i+1}(e_{i-1}e_i-q_{i}^{-1}e_ie_{i-1})e_{i+1}  - (e_{i-1}e_i-q_{i}^{-1}e_ie_{i-1})e_{i+1}e_{i+2}e_{i+1} \\
&\quad + {(-1)^{\e_{i+1}}[2]}e_{i+1}(e_{i-1}e_i-q_{i}^{-1}e_ie_{i-1})e_{i+2}e_{i+1},
\end{split}
\end{equation*}
which is zero by \eqref{eq:Serre-rel-2} for $\U(\e)$ with respect to $i+1 \in I_{\rm odd}$.
The proof for $(i-3,i-2,i-1)$ is the same. The proof for $(i-2,i-1,i)$ reduces to the previous cases by applying $\tau_i$ to $E_l$ for $l=i-2,i-1,i$. We leave the proof for $F_j$'s to the reader.
\qed
\vskip 2mm


\subsection{Truncation to $\U(\e')$-modules}
Let $\e'$ be as in Section \ref{subsec:subalg e'}. Suppose that $M'$ is the number of $j$'s with $\e'_j=0$ and $N'$ is the number of $j$'s with $\e'_j=1$ in $\e'$.

 For a submodule $V$ of $\mc{V}^{\ot \ell}$ ($\ell\geq1$), we define
\begin{equation}\label{eq:truncation-1}
\mf{tr}^\e_{\e'}(V) = \bigoplus_{\substack{\mu\in {\rm wt}(V) \\ (\mu|\de_i)=0}}V_\mu,
\end{equation}
where ${\rm wt}(V)$ is the set of weights of $V$. For any submodules $V,W$ of a tensor power of $\V$, it is clear that 
\begin{equation*}
\mf{tr}^\e_{\e'}(V\ot W)=\mf{tr}^\e_{\e'}(V)\ot \mf{tr}^\e_{\e'}(W),
\end{equation*}
as a vector space.

\begin{lem}\label{lem:truncation of natural repn}
Let $\V'=\mf{tr}^\e_{\e'}(\V)$. Then
\begin{itemize}
\item[(1)] $\V'$ is isomorphic to the natural representation of $\ring\U(\e')$ given in \eqref{eq:natural repn},

\item[(2)] $\mf{tr}^\e_{\e'}(\V^{\ot \ell})$ is isomorphic to $\V'^{\ot \ell}$ as a $\ring\U(\e')$-module.

\end{itemize}
\end{lem}
\pf (1) Let us assume that $2\leq i\leq n-2$ since the proof for the other cases is similar.
Let $j\in \ring {(I')}$ and $k\in \I\setminus\{i\}$ given. It is clear from \eqref{eq:generators for e'} that
\begin{equation*}\label{eq:natural for e'-1}
E_jv_k=
\begin{cases}
v_j, & \text{if $k= j+1$},\\
0, & \text{if $k\neq j+1$},
\end{cases} \quad (j\leq i-2),\quad
E_jv_k=
\begin{cases}
v_{{j+1}}, & \text{if $k= j+2$},\\
0, & \text{if $k\neq j+2$},
\end{cases}\quad (j\geq i).
\end{equation*}
When $j=i-1$, we have $E_{i-1}=e_{i-1}e_i-q^{-1}_ie_ie_{i-1}$, and 
\begin{equation*}\label{eq:natural for e'-2}
E_{i-1}v_k=
\begin{cases}
v_{i-1}, & \text{if $k= i+1$},\\
0, & \text{if $k\neq i+1$}.
\end{cases}
\end{equation*}
We have similar formulas for $F_j$ for $j\in \ring {(I')}$. Hence $\V'$ is invariant under the action of $\ring\U(\e')$. In fact, $\mc{V}'$ is isomorphic to the natural representation of $\ring\U(\e')$ \eqref{eq:natural repn}.

(2) We see that the actions of $E_j, F_j, K_j$ $(j\in \ring {(I')})$ on $\V'\ot\V'$ are equal to those of 
\begin{equation}\label{eq:action of e' on tensor power}
K_j^{-1} \ot E_j + E_j\ot 1 ,\quad 1\ot F_j + F_j\ot K_j, \quad K_j\ot K_j, 
\end{equation}
respectively. 
This implies that $\V'\ot\V'$ and hence $(\V')^{\ot \ell}$ are invariant under the action of $\ring\U(\e')$.
For example, in case of $E_{i-1}=e_{i-1}e_i-q^{-1}_ie_ie_{i-1}$, we have
\begin{equation*}
\begin{split}
\Delta(E_{i-1})
&=\Delta(e_{i-1})\Delta(e_{i})-q^{-1}_i\Delta(e_{i})\Delta(e_{i-1})\\
&=k_{i-1}^{-1}k_i^{-1}\ot e_{i-1}e_i + k_{i-1}^{-1}e_i\ot e_{i-1}+ k_{i}^{-1}e_{i-1}\ot e_{i} + e_{i-1}e_i\ot 1 \\
&\ \ -q^{-1}_i 
( 
k_i^{-1}k_{i-1}^{-1}\ot e_ie_{i-1} + k_{i}^{-1}e_{i-1}\ot e_{i}+ k_{i-1}^{-1}e_{i}\ot e_{i-1} + e_{i}e_{i-1}\ot 1
).
\end{split}
\end{equation*}
Then the action of $\Delta(E_{i-1})$ on $\mc{V}'\ot \mc{V}'$ is equal to $k_{i}^{-1}k_{i-1}^{-1}\ot e_{i-1}e_i + e_{i-1}e_{i}\ot 1$, and hence $K_{i-1}^{-1}\ot E_{i-1} + E_{i-1}\ot 1$.
\qed

\begin{prop}\label{prop:truncation of poly}
Let $\lambda \in \cP_{M|N}$ be given.
\begin{itemize}
\item[(1)] $\mf{tr}^\e_{\e'}(V_\e(\la))$ is a $\ring\U(\e')$-submodule of $V_\e(\la)$ via $\phi$.

\item[(2)] $\mf{tr}^\e_{\e'}(V_\e(\la))$ is non-zero if and only if $\la\in \cP_{M'|N'}$. 
In this case, we have 
\begin{equation*}
\mf{tr}^\e_{\e'}(V_\e(\la)) \cong V_{\e'}(\la),
\end{equation*}
as a $\ring\U(\e')$-module.  

\end{itemize}
\end{prop}
\pf (1) It follows from Lemma \ref{lem:truncation of natural repn} and
\begin{equation}\label{eq:truncation-2}
\mf{tr}^\e_{\e'}(V_\e(\la))=\mf{tr}^\e_{\e'}(Y^\la(q){\mc V}^{\ot \ell}) =Y^\la(q)\mf{tr}^\e_{\e'}({\mc V}^{\ot \ell})
=Y^\la(q)\mf{tr}^\e_{\e'}({\mc V})^{\ot \ell}=Y^\la(q)(\V')^{\ot \ell}.
\end{equation}

(2) Note that $SST_{\e'}(\la)\subset SST_{\e}(\la)$.
By Proposition \ref{prop:poly repn} and \eqref{eq:truncation-2}, we see that $\mf{tr}^\e_{\e'}(V_\e(\la))$ is a $\Q(q)$-span of $\{\,v_{T'}\,|\, T'\in SST_{\e'}(\la)\,\}$, which in fact forms a basis. This implies that $\mf{tr}^\e_{\e'}(V_\e(\la))$ is non-zero if and only if $\la\in \cP_{M-1|N}$ when $\e_i=0$, and $\la\in \cP_{M|N-1}$ when $\e_i=1$.
Hence, $\mf{tr}^\e_{\e'}(V_\e(\la))$ is isomorphic to $V_{\e'}(\la)$ when it is non-zero by  \eqref{eq:action of e' on tensor power} and Proposition \ref{prop:poly repn}. \qed

\begin{cor}
Let $V,W$ be submodules of a tensor power of $\V$. Then 
\begin{itemize}
\item[(1)] $\mf{tr}^\e_{\e'}(V)$, $\mf{tr}^\e_{\e'}(W)$, and $\mf{tr}^\e_{\e'}(V\ot W)$ are $\ring\U(\e')$-modules via $\phi$,

\item[(2)] $\mf{tr}^\e_{\e'}(V\ot W)\cong\mf{tr}^\e_{\e'}(V)\ot \mf{tr}^\e_{\e'}(W)$ as $\ring\U(\e')$-modules. 
\end{itemize}
\end{cor}
\pf Since $\V^{\ot \ell}$ is completely reducible, it follows from Proposition \ref{prop:truncation of poly} and \eqref{eq:action of e' on tensor power}.
\qed
\vskip 2mm

We may define $\mf{tr}^\e_{\e'}$ and have similar results for $\U(\e)$-modules in $\mc O_{\geq 0}$.

\begin{prop}\label{prop:truncation of fund}
\mbox{} 
\begin{itemize}
\item[(1)] For $s\in\Z_+$ and $x\in \Q(q)$, $\mf{tr}^\e_{\e'}(\W_{s,\e}(x))$ is a $\U(\e')$-submodule of $\W_{s,\e}(x)$ via $\phi$, and 
\begin{equation*}
\mf{tr}^\e_{\e'}(\W_{s,\e}(x)) \cong \W_{s,\e'}(x).
\end{equation*}
Moreover, $(\mf{tr}^\e_{\e'}(\mc{L}_{s,\e}),\mf{tr}^\e_{\e'}(\mc{B}_{s,\e}))$ is a crystal base of $\mf{tr}^\e_{\e'}(\W_{s,\e})$ isomorphic to $(\mc{L}_{s,\e'},\mc{B}_{s,\e'})$.
\item[(2)] For $l,m\in\Z_+$ and $x,y\in \Q(q)$, $\mf{tr}^\e_{\e'}(\W_{l,\e}(x)\ot \W_{m,\e}(y))$ is a $\U(\e')$-module via $\phi$, and
\begin{equation*}
\mf{tr}^\e_{\e'}(\W_{l,\e}(x)\ot \W_{m,\e}(y))\cong\mf{tr}^\e_{\e'}(\W_{l,\e}(x))\ot \mf{tr}^\e_{\e'}(\W_{m,\e}(y)),
\end{equation*}
as $\U(\e')$-modules.
\end{itemize}
\end{prop}
\pf The proof is the same as in Proposition \ref{prop:truncation of poly}.
\qed

\subsection{Irreducibility of $\W_{l,\e}(x)\ot \W_{m,\e}(y)$}
Let us show that $\W_{l,\e}(x)\ot \W_{m,\e}(y)$ is irreducible for $l,m\in\Z_+$ and generic $x,y\in \Q(q)$.
When $\e=\e_{M|N}$, the irreducibility is shown in \cite{KOS}. In this paper, we give a different proof of it, which is also available for arbitrary $\e$.

\begin{thm}\label{thm:main theorem-1}
For $l,m\in\Z_+$, $\W_{l,\e}\ot \W_{m,\e}$ is irreducible.
\end{thm}
\pf Let us assume without loss of generality that $M, N\geq 1$ with $\e_1=0$.

Let $(\mc{L}_{s,\e},\mc{B}_{s,\e})$ be the crystal base of $\W_{s,\e}$ in \eqref{eq:crystal base of W(r)} for $s=l,m$. By Proposition \ref{prop:tensor product rule}, $(\mc{L}_{l,\e}\ot \mc{L}_{m,\e},\mc{B}_{l,\e}\ot \mc{B}_{m,\e})$ is a crystal base of $\W_{l,\e}\ot \W_{m,\e}$.
If $M=1$, then it is proved in \cite{KO} that $\mc{B}_{l,\e}\ot \mc{B}_{m,\e}/\{\pm 1\}$ is connected.  Since $\dim_{\Q(q)}(\W_{l,\e}\ot \W_{m,\e})_{(l+m)\de_1}=1$ and $\mc{B}_{l,\e}\ot \mc{B}_{m,\e}/\{\pm 1\}$ is connected, it follows from \cite[Lemma 2.7]{BKK} that $\W_{l,\e}\ot \W_{m,\e}$ is irreducible.

We assume that $M\geq 2$. Set $\e'=\e_{M|0}$, which is the subsequence of $\e$ obtained by removing all $\e_i=1$'s. Note that the length of $\e'$ may be less than 4 so that $\U(\e')$ is not well-defined, but $\mf{tr}^\e_{\e'}$ can be defined in the same way as in \eqref{eq:truncation-1}. We put
\begin{equation*}
\mc{W}_{s,\e'}:=\mf{tr}^\e_{\e'}(\mc{W}_{s,\e}),\quad 
\mc{L}_{s,\e'}:=\mf{tr}^\e_{\e'}(\mc{L}_{s,\e})\subset \mc{L}_{s,\e},\quad 
\mc{B}_{s,\e'}:=\mf{tr}^\e_{\e'}(\mc{B}_{s,\e})\subset \mc{B}_{s,\e}.
\end{equation*}

Let $1\leq j_1<\dots<j_M\leq n$ be such that $\e_{j_k}=0$ for $1\leq k\leq M$. 
By Theorem \ref{thm:folding homomorphism}, we have a $U_q(\mf{sl}_2)$-action on $\W_{l,\e'}\ot \W_{m,\e'}$ corresponding to the pair $(\e_{j_k},\e_{j_{k+1}})$ or $(\e_{j_M},\e_{j_1})$. 
For $0\leq k\leq M-1$, let us denote by $\te_{k'}$ and $\tf_{k'}$ the Kashiwara operators corresponding to $(\e_{j_k},\e_{j_{k+1}})$ when $k\neq 0$ and to $(\e_{j_M},\e_{j_1})$ when $k=0$. 

If we put $I'=\{\,k'\,|\,k=0,\dots,M-1 \,\}$, then $\mc{L}_{l,\e'}\ot \mc{L}_{m,\e'}$ is invariant under $\te_{k'}$ and $\tf_{k'}$ for $k'\in I'$, and hence $\mc{B}_{l,\e'}\ot \mc{B}_{m,\e'}/\{\pm 1\}$ is an $I'$-colored oriented graph.
Since $\mc{L}_{l,\e'}\ot \mc{L}_{m,\e'}\subset \mc{L}_{l,\e}\ot \mc{L}_{m,\e}$ and $\mc{B}_{l,\e'}\ot \mc{B}_{m,\e'}\subset \mc{B}_{l,\e}\ot \mc{B}_{m,\e}$, we may regard $\mc{B}_{l,\e}\ot \mc{B}_{m,\e}/\{\pm 1\}$ as an $(I\sqcup I')$-colored oriented graph.

Let ${\bf b}=|{\bf m}_1\rangle\ot |{\bf m}_2\rangle\in \mc{B}_{l,\e}\ot \mc{B}_{m,\e}$ be given. We will show that ${\bf b}$ is connected to $|l\be_1\rangle\ot |m\be_1\rangle$, which implies that  $\mc{B}_{l,\e}\ot \mc{B}_{m,\e}/\{\pm 1\}$ is connected as an $(I\sqcup I')$-colored oriented graph.
Let us write ${\bf m}_i=(m_{i1},\dots,m_{in})$ for $i=1,2$.

We first claim that there exists a sequence $i_1,\dots,i_r\in I$ such that $(\e_{i_k},\e_{i_{k}+1})\neq (0,0)$ for $1\leq k\leq r$ and
\begin{equation}\label{eq:path to usual affine type}
{\bf b}':=\td{x}_{i_1}\dots\td{x}_{i_r} 
{\bf b} \equiv 
|{\bf m}'_1\rangle\ot |{\bf m}'_2\rangle \pmod{q\mc{L}_{l,\e'}\ot \mc{L}_{m,\e'}},
\end{equation}
for some $|{\bf m}'_1\rangle\in \W_{l,\e'}$ and $|{\bf m}'_2\rangle\in \W_{m,\e'}$, where $\td{x}_{i_s}=\te_{i_s}$ or $\tf_{i_s}$ for each $1\leq s\leq r$.

Suppose that there exists $k$ with $\e_k=1$ such that $m_{1k}=1$ or $m_{2k}=1$. 
Let $i$ and $j$ be the maximal and minimal indices respectively such that $i<k<j$ and $\e_i=\e_j=0$. If there is no such $(i,j)$, then we have $\e=\e_{M|N}$ and identify this case with the one of $\e=(0^{M-1},1^N,0)$. 
Since we will choose $i_1,\dots, i_r$ in $\{i,i+1,\dots,j-1\}$, we may assume for simplicity that $m_{ab}=0$ for $a=1,2$ and $b\not\in \{i,\dots, j\}$.

Let us use induction on $L=|{\bf m}_1|+|{\bf m}_2|$. 
Suppose that $L=1$. 
If $m_{1k}=1$, then $\tf_{j-1}\tf_{j-2}\dots\tf_{k}{\bf b}$ satisfies \eqref{eq:path to usual affine type}. 
If $m_{2k}=1$, then $\te_{i}\te_{i+1}\dots\te_{k-1}{\bf b}$ satisfies \eqref{eq:path to usual affine type}. 

Suppose that $L>1$.
We may assume that $\tf_{i+1} {\bf b}=\tf_{i+2} {\bf b}= \dots =\tf_{j-1} {\bf b}=0$. 
Then by tensor product rule in Proposition \ref{prop:tensor product rule} we have
\begin{equation}\label{eq:m step1}
\begin{split}
{\bf m}_1&= m_{1i}\be_i + \sum_{x\leq u\leq y}\be_u + \sum_{z\leq v\leq j-1}\be_v + m_{1j}\be_j , \\
{\bf m}_2&= m_{2i}\be_i + \sum_{y+1\leq v\leq j-1}\be_v + m_{2j}\be_j,
\end{split}
\end{equation}
for some $i<x<y<z<j$. Here we assume that $\sum_{z\leq v\leq j-1}\be_v$ in ${\bf m}_1$ is empty if there is no such $z$. Now we take the following steps to construct ${\bf b}'$ in \eqref{eq:path to usual affine type}.

{Step 1}. If there exists $z$ such that $y<z<j$ and $m_{1z}=\dots=m_{1j-1}=1$, then by applying $\tf_{z}\tf_{z-1}\dots\tf_{j-1}$ to ${\bf b}$, ${\bf m}_1$ in \eqref{eq:m step1} is replaced by 
\begin{equation}\label{eq:m step1-2}
m_{1i}\be_i + \sum_{x\leq u\leq y}\be_u + \sum_{z+1\leq v\leq j-1}\be_v + (m_{1j}+1)\be_j.
\end{equation}
Repeating this step, \eqref{eq:m step1-2} is replaced by
\begin{equation*}\label{eq:m step 1-3}
m_{1i}\be_i + \sum_{x\leq u\leq y}\be_u + (m_{1j}+j-z)\be_j.
\end{equation*}
Hence we may assume that ${\bf m}_1$ in \eqref{eq:m step1} is of the form $m_{1i}\be_i + \sum_{x\leq u\leq y+1}\be_u + m_{1j}\be_j$.

{Step 2}.  If $m_{1j}=0$, then we have
\begin{equation*}
\tf_{j-1}{\bf b}=|{\bf m}_1\rangle \ot |{\bf m}_2-\be_{j-1}+\be_{j}\rangle.
\end{equation*}
Hence we may apply the induction hypothesis to conclude \eqref{eq:path to usual affine type}.

{Step 3}. If $m_{ij}\neq 0$, then by applying $\te_i\te_{i+1}\dots\te_{j-2}\te_{j-1}$ to ${\bf b}$, ${\bf m}_1$ and ${\bf m}_2$ are replaced by
\begin{equation*}\label{eq:m step 3-1}
\begin{split}
& m_{1i}\be_i + \sum_{x\leq u\leq y+1}\be_u + (m_{1j}-1)\be_j, \\
& (m_{2i}+1)\be_i + \sum_{y+2\leq v\leq j-1}\be_v + m_{2j}\be_j,
\end{split}
\end{equation*}
respectively. Repeating this step $d$ times such that $m_{1j}-d\geq 0$ and $y+d+1\leq j$, ${\bf m}_1$ and ${\bf m}_2$ are replaced by
\begin{equation*}\label{eq:m step 3-2}
\begin{split}
& m_{1i}\be_i + \sum_{x\leq u\leq y+d}\be_u + (m_{1j}-d)\be_j, \\
& (m_{2i}+d)\be_i + \sum_{y+d+1\leq v\leq j-1}\be_v + m_{2j}\be_j,
\end{split}
\end{equation*}
respectively. 
We may keep this process until $m_{1j}-d=0$, which belongs to the case in {\em Step 2}, or $\sum_{y+d+1\leq v\leq j-1}\be_v$ is empty. 
In the latter case, ${\bf m}_1$ is replaced by $m_{1i}\be_i + \sum_{x\leq u\leq j-1}\be_u + (m_{1j}-d)\be_j$ so that we may apply $\tf_{j-1}$ and use induction hypothesis to have ${\bf b}'$.
This proves the claim.

By construction of ${\bf b}'$ and its weight, we have
\begin{equation*}
{\bf b}'- |{\bf m}'_1\rangle\ot |{\bf m}'_2\rangle\in \left(\mc{L}_{l,\e'}\ot \mc{L}_{m,\e'}\right)\cap \left(q\mc{L}_{l,\e}\ot \mc{L}_{m,\e}\right)=q\mc{L}_{l,\e'}\ot \mc{L}_{m,\e'},
\end{equation*}
and hence ${\bf b}'\in (\mc{L}_{l,\e'}\ot \mc{L}_{m,\e'}/q\mc{L}_{l,\e'}\ot \mc{L}_{m,\e'})\subset (\mc{L}_{l,\e}\ot \mc{L}_{m,\e}/q\mc{L}_{l,\e}\ot \mc{L}_{m,\e})$. 
If $M=2$, then it is easy to show that ${\bf b}'=|{\bf m}'_1\rangle\ot |{\bf m}'_2\rangle\in \mc{B}_{l,\e'}\ot \mc{B}_{m,\e'}$ is connected to $|l\be_1\rangle\ot |m\be_1\rangle$ under $\te_{k'}$ and $\tf_{k'}$ for $k=0,1$.
If $M\geq 3$, then we can also show that ${\bf b}'=|{\bf m}'_1\rangle\ot |{\bf m}'_2\rangle\in \mc{B}_{l,\e'}\ot \mc{B}_{m,\e'}$ is connected to $|l\be_1\rangle\ot |m\be_1\rangle$ by using the fact that  $\mc{B}_{l,\e'}\ot \mc{B}_{m,\e'}/\{\pm 1\}$ is a connected crystal of type $A_{M-1}^{(1)}$ (cf.~\cite{AK}).

Finally, since $\dim_{\Q(q)}(\W_{l,\e}\ot \W_{m,\e})_{(l+m)\de_1}=1$ and $\mc{B}_{l,\e}\ot \mc{B}_{m,\e}/\{\pm 1\}$ is connected, it follows from \cite[Lemma 2.7]{BKK} that $\W_{l,\e}\ot \W_{m,\e}$ is irreducible.
This completes the proof. \qed

\begin{cor}\label{cor:irreducibility}
For $l,m\in\Z_+$ and generic $x,y\in \Q(q)$, $\W_{l,\e}(x)\ot \W_{m,\e}(y)$ is irreducible.
\end{cor}
\pf It follows from \cite[Lemma 3.4.2]{KMN}.
\qed

\subsection{Existence of $R$ matrix}
For $l,m\in\Z_+$ and $x,y\in\Q(q)$,
consider a non-zero $\Q(q)$-linear map $R$ on $\W_{l,\e}(x)\otimes \W_{m,\e}(y)$ such that
\begin{equation}\label{eq:R matrix}
\begin{split}
\Delta^{\rm op}(g) \circ R = R \circ  \Delta (g),
\end{split}
\end{equation}
for $g\in \U(\e)$, where $\Delta^{\rm op}$ is the opposite coproduct of $\Delta$ in \eqref{eq:comult-1}, that is, $\Delta^{\rm op}(g) = P \circ \Delta(g) \circ P $ and $P(a\otimes b)=b\otimes a$. We denote it by $R(z)$, where $z=x/y$, since $R$ depends only on $z$. 

We say that $R(z)$ satisfies the Yang-Baxter equation if we have
\begin{equation}\label{eq:YB equation}
R_{12}(u)R_{13}(uv)R_{23}(v)=R_{23}(v)R_{13}(uv)R_{12}(u),
\end{equation}
on $\W_{s_1,\e}(x_1)\otimes \W_{s_2,\e}(x_2)\otimes \W_{s_3,\e}(x_3)$ with $u=x_1/x_2$ and $v=x_2/x_3$ for $(s_1),(s_2),(s_3)\in \cP_{M|N}$. Here $R_{ij}(z)$ denotes the map which acts as $R(z)$ on the $i$-th and the $j$-th component and the identity elsewhere. We call $R(z)$ the {\em (quantum) $R$ matrix}.

\begin{thm}\label{thm:existence of R-matrix}
Let $l, m\in \Z_+$ given with $(l), (m)\in \cP_{M|N}$. 
Suppose that $\e_1=0$.
There exists a unique non-zero linear map $R(z)\in {\rm End}_{\Q(q)}(\W_{l,\e}(x)\otimes \W_{m,\e}(y))$ satisfying \eqref{eq:R matrix} and \eqref{eq:YB equation}, and $R(z)(|l\be_1\rangle\ot |m\be_1\rangle)=|l\be_1\rangle\ot |m\be_1\rangle$ for generic $x,y\in \Q(q)$.
\end{thm}
\pf The existence of such a map for arbitrary $\e$ is proved in \cite[Theorem 5.1]{KOS} with respect to $\Delta_+$ in  \eqref{eq:comult-2}, say $R_+$. 
Let 
\begin{equation}\label{eq:chi}
\chi=\psi\circ(\phi\ot \phi),
\end{equation}
where $\psi$ and $\phi$ are given in \eqref{eq:isomorphism for two comult} and \eqref{eq:iso phi}, respectively.
Then
\begin{equation*}
R:=\chi^{-1}\circ R_+\circ \chi
\end{equation*}
satisfies the conditions \eqref{eq:R matrix} and \eqref{eq:YB equation}, and $R(z)(|l\be_1\rangle\ot |m\be_1\rangle)=|l\be_1\rangle\ot |m\be_1\rangle$ with respect to $\Delta$. The uniqueness follows from the irreducibility in Corollary \ref{cor:irreducibility} and normalization by $R(z)(|l\be_1\rangle\ot |m\be_1\rangle)=|l\be_1\rangle\ot |m\be_1\rangle$.
\qed

\begin{rem}{\rm
If $M=0$, then the existence of $R$ matrix is already known. Hence we may assume that $M\geq 1$. If $\e_1=1$, then we may choose the smallest $i\in \I_0$ so that there exists a unique $R$ matrix satisfying $R(z)(|l\be_i\rangle\ot |m\be_i\rangle)=|l\be_i\rangle\ot |m\be_i\rangle$.

}
\end{rem}

\section{Kirillov-Reshetikhin modules}\label{sec:KR modules}

\subsection{Spectral decomposition}
Suppose that $\e=(\e_1,\dots,\e_n)$ is given with $n\geq 4$. 
Recall that $M$ is the number of $i$'s with $\e_i=0$ and $N$ is the number of $i$'s with $\e_i=1$ in $\e$.

Let $l,m\in\Z_+$ be given. 
Let $R_\e(z)$ be the $R$ matrix on $\W_{l,\e}(x)\otimes \W_{m,\e}(y)$ in Theorem \ref{thm:existence of R-matrix}.
We have as a $\ring\U(\e)$-module,
\begin{equation*}\label{eq:decomp of W_l and W_m}
\W_{l,\e}(x)\otimes \W_{m,\e}(y) \cong 
\bigoplus_{t\in H(l,m)} V_\e((l+m-t,t)),
\end{equation*}
where $H(l,m)=\{\,t\,|\,0\leq t\leq \min\{l,m\}, (l+m-t,t)\in \cP_{M|N}\,\}$.

Let us take a sequence $\e''=(\e''_1,\dots,\e''_{n''})$ of $0,1$'s with $n''\gg n$ satisfying the following:
\begin{itemize}
\item[(1)] $\e$ is a subsequence of $\e''$, 

\item[(2)] we have as a $\ring\U(\e'')$-module
\begin{equation*}\label{eq:decomp of W_l and W_m-e''}
\W_{l,\e''}(x)\otimes \W_{m,\e''}(y) \cong 
\bigoplus_{0\leq t\leq \min\{l,m\}} V_{\e''}((l+m-t,t)),
\end{equation*}

\item[(3)] if $\e'=\e_{M''|0}$ with $M''=|\{\,i\,|\,\e''_i=0\,\}|\gg 0$, then we have as a $\ring\U(\e')$-module
\begin{equation*}\label{eq:decomp of W_l and W_m-e'}
\W_{l,\e'}(x)\otimes \W_{m,\e'}(y) \cong 
\bigoplus_{0\leq t\leq \min\{l,m\}} V_{\e'}((l+m-t,t)),
\end{equation*}
\end{itemize}

Let $R_{\e''}(z)$ and $R_{\e'}(z)$ denote the $R$ matrices on 
$\W_{l,\e''}(x)\otimes \W_{m,\e''}(y)$ and $\W_{l,\e'}(x)\otimes \W_{m,\e'}(y)$, respectively.

\begin{lem}\label{lem:comm diagram of R matrix}
For $\upepsilon=\e$ or $\e'$, we have the following commutative diagram:
\begin{equation*}\label{eq:comm diagram of R matrix}
\xymatrixcolsep{4pc}\xymatrixrowsep{2pc}\xymatrix{
\W_{l,\e''}(x)\otimes \W_{m,\e''}(y) \ar@{->}^{PR_{\e''}(z)}[r] \ar@{->}^{\mf{tr}^{\e''}_{\upepsilon}}[d] &\ \W_{m,\e''}(y)\otimes \W_{l,\e''}(x) \ar@{->}^{\mf{tr}^{\e''}_{\upepsilon}}[d] \\  
\W_{l,\upepsilon}(x)\otimes \W_{m,\upepsilon}(y) \ar@{->}^{PR_{\upepsilon}(z)}[r] &\  \W_{m,\upepsilon}(y)\otimes \W_{l,\upepsilon}(x)}\quad 
\end{equation*}
\end{lem}
\pf
For $\upepsilon=\e$ or $\e'$, the restriction of $PR_{\e''}(z)$ on $\mf{tr}^{\e''}_{\upepsilon}\left(\W_{l,\e}(x)\otimes \W_{m,\e}(y)\right)$, which gives a well-defined $\U(\upepsilon)$-linear endomorphism. By Proposition \ref{prop:truncation of fund} and Theorem \ref{thm:existence of R-matrix}, the restricted $R$ matrix is the quantum  $R$ matrix on $\W_{l,\upepsilon}(x)\otimes \W_{m,\upepsilon}(y)$, which proves the commutativity of the diagram.
\qed

For $0\leq t\leq \min\{l,m\}$, let $v'(l,m,t)$ be the highest weight vectors of $V_{\e'}((l+m-t,t))$ in $\W_{l,\e'}(x)\otimes \W_{m,\e'}(y)$ such that 
\begin{equation*}
\begin{split}
& v'(l,m,t)\in \mc{L}_{l,\e'}\ot \mc{L}_{m,\e'},\\
& v'(l,m,t)\equiv | l\be_1 \rangle\ot | (m-t)\be_1 + t\be_2 \rangle \pmod{q\mc{L}_{l,\e'}\ot \mc{L}_{m,\e'}}.
\end{split}
\end{equation*}
We also define $v'(m,l,t)$ in the same manner. 
For $0\leq t'\leq \min\{l,m\}$, we may regard 
\begin{equation*}
V_{\e}((l+m-t,t))\subset V_{\e''}((l+m-t,t)), \quad V_{\e'}((l+m-t,t))\subset V_{\e''}((l+m-t,t))
\end{equation*}
as a $\Q(q)$-space, and let $\mc P^{l,m}_{t} : \W_{l,\e''}(x)\otimes \W_{m,\e''}(y) \longrightarrow \W_{m,\e''}(y)\otimes \W_{l,\e''}(x)$ be a $\ring\U(\e'')$-linear map given by $\mc P^{l,m}_{t}(v'(l,m,t'))=\delta_{tt'}v'(m,l,t')$. Then we have the following spectral decomposition of $PR_{\e''}(z)$
\begin{equation*}
PR_{\e''}(z)= \sum_{0\leq t\leq \min\{l,m\}}\rho_t(z) \mc P^{l,m}_t,
\end{equation*}
for some $\rho_t(z)\in \mathbb{Q}(q)$. 
By Proposition \ref{prop:truncation of poly} and Lemma \ref{lem:comm diagram of R matrix}, we also have the following spectral decomposition of $PR_{\upepsilon}(z)$
\begin{equation}\label{eq:R for e and e'}
\begin{split}
PR_{\e'}(z)&= \sum_{0\leq t\leq \min\{l,m\}}\rho_t(z) \mc P^{l,m}_t,\\
PR_{\e}(z)&= \sum_{t\in H(l,m)}\rho_t(z) \mc P^{l,m}_t,
\end{split}
\end{equation}
where we understand $\mc{P}^{l,m}_t$ as defined on $\W_{l,\upepsilon}(x)\otimes \W_{m,\upepsilon}(y)$. 
Then we have the following explicit description of $PR_\e(z)$, which is proved in case of $\e=\e_{M|N}$ \cite{KOS}.

\begin{thm}\label{thm:spectral decomp}
We have

\begin{equation}\label{eq:spectral M=0}
PR_\e(z) =  
\sum_{t=\max\{l+m-n,0\}}^{\min\{l,m\}}\left(\prod_{i=t+1}^{\min\{l,m\}}\frac{z-q^{l+m-2i+2}}{1-q^{l+m-2i+2}z}\right)\mc P^{l,m}_t \qquad (M=0),
\end{equation}

\begin{equation}\label{eq:spectral M=1}
PR_\e(z)=
\sum_{t=0}^{\min\{l,m,n-1\}}\left(\prod_{i=1}^{t}\frac{1-q^{l+m-2i+2}z}{z-q^{l+m-2i+2}}\right)\mc P^{l,m}_t\hskip 2cm (M=1),
\end{equation}

\begin{equation}\label{eq:spectral M>1}
PR_\e(z)=
\sum_{t=0}^{\min\{l,m\}}
\left(\prod_{i=1}^{t}\frac{1-q^{l+m-2i+2}z}{z-q^{l+m-2i+2}}\right)\mc P^{l,m}_t
\hskip 2cm (2\leq M\leq n),
\end{equation}
where we assume that $\rho_{\min\{l,m\}}(z)=1$ in \eqref{eq:spectral M=0} and $\rho_0(z)=1$ in \eqref{eq:spectral M=1} and \eqref{eq:spectral M>1}.
\end{thm}
\pf We may consider the case of $1\leq M\leq n$ only since the case when $M=0$ is known (see \cite[(5.6)]{KO} or \cite[(6.10)]{KOS}). 
It is well-known that $PR_{\e'}(z)$ for $\e'=\e_{M''|0}$ has the following spectral decomposition
\begin{equation*}
PR_{\e'}(z)= \sum_{0\leq t\leq \min\{l,m\}}\rho'_t(z) \mc P^{l,m}_t,
\end{equation*}
where
\begin{equation*}
\rho'_0(z)=1,\quad \rho'_t(z)=\prod_{i=1}^{t}\frac{1-q^{l+m-2i+2}z}{z-q^{l+m-2i+2}} \qquad(1\leq t\leq \min\{l,m\}),
\end{equation*}
(cf.~\cite[(5.8)]{KO} or \cite[(6.16)]{KOS}). We remark that $\chi(v'(l,m,t))$ and $\chi(v'(m,l.t))$ for $0\leq t\leq \min\{l,m\}$ are the same scalar multiplications of the highest weight vectors in \cite[(6.14)]{KOS}, where $\chi$ is as in \eqref{eq:chi}.
Hence it follows from \eqref{eq:R for e and e'} that
\begin{equation*}
\rho_t(z)=\rho'_t(z) \quad (t\in H(l,m)),
\end{equation*}
which completes the proof.
\qed

\subsection{Kirillov-Reshetikhin modules}

As an application of Theorem \ref{thm:spectral decomp}, let us construct a family of irreducible $\U(\e)$-modules in $\mc{O}_{\geq 0}$ which corresponds to usual Kirillov-Reshetikhin modules under truncation.
Let us assume that $1\leq M \leq n-1$ since the results when $M\in \{0, n\}$ are well-known \cite{KMN2}. 

Fix $s\geq 1$ and put $V_x=\mc{W}_{s,\e}(x)$ for $x\in \Q(q)$.
We take a normalization
\begin{equation*}
\check{R}(z) = \left(\prod_{i=1}^s \frac{z-q^{2s-2i+2}}{1-q^{2s-2i+2}z}\right) PR(z),
\end{equation*}
where $R(z)$ is the $R$ matrix on $V_x\otimes V_y$.
Since $(s^2)\not\in \cP_{M|N}$ if and only if $M=1$ and $s>n-1$, we have
\begin{equation*}\label{eq:normalized R}
\check{R}(z)=
\begin{cases}
\sum_{t=0}^{n-1}\left(\prod_{i=t+1}^{s}\dfrac{z-q^{2s-2i+2}}{1-q^{2s-2i+2}z}\right)\mc P^{s,s}_t,& \text{if $(s^2)\not\in \cP_{M|N}$},\\
\mc P^{s,s}_s + \sum_{t=0}^{s-1}\left(\prod_{i=t+1}^{s}\dfrac{z-q^{2s-2i+2}}{1-q^{2s-2i+2}z}\right)\mc P^{s,s}_t,& \text{if $(s^2)\in \cP_{M|N}$}. 
\end{cases}
\end{equation*}

For $r\geq 2$, let $W$ denote the group of permutations on $r$ letters generated by $s_i=(i\ i+1)$ for $1\leq i\leq r-1$.
By Theorem \ref{thm:existence of R-matrix}, we have $\U(\e)$-linear maps 
\begin{equation*}
\check{R}_w(x_1,\ldots,x_r) : 
V_{x_1}\otimes \cdots\otimes V_{x_r} \longrightarrow V_{x_{w(1)}}\otimes \cdots\otimes V_{x_{w(r)}}
\end{equation*}
for $w\in W$ and generic $x_1,\ldots, x_r$ satisfying the following:
\begin{equation*}\label{eq:R matrix for w}
\begin{split}
&\check{R}_1(x_1,\ldots,x_r) = {\rm id}_{V_{x_1}\otimes \cdots\otimes V_{x_r}},\\
&\check{R}_{s_i}(x_1,\ldots,x_r) = \left(\otimes_{j<i}{\rm id}_{V_{x_j}}\right)\otimes \check{R}(x_i/x_{i+1}) \otimes \left(\otimes_{j>i+1}{\rm id}_{V_{x_j}}\right),\\
&\check{R}_{ww'}(x_1,\ldots,x_r) = \check{R}_{w'}(x_{w(1)},\ldots,x_{w(r)})\check{R}_{w}(x_1,\ldots,x_r),
\end{split}
\end{equation*}
for $w, w'\in W$ with $\ell(ww')=\ell(w)+\ell(w')$.
Let $w_0$ denote the longest element in $W$. 
By Theorem \ref{thm:spectral decomp}, $\check{R}_{w_0}(x_1,\ldots,x_r)$ does not have a pole at $q^{2k}$ for $k\in \Z_+$ as a function in $x_1,\ldots,x_r$. Hence we have a $\U(\e)$-linear map
\begin{equation*}
\begin{split}
\check{R}_r := \check{R}_{w_0}(q^{r-1},q^{r-3},\cdots, q^{1-r}) : V_{q^{r-1}}\otimes \cdots\otimes V_{q^{1-r}} \longrightarrow 
V_{q^{1-r}}\otimes \cdots\otimes V_{q^{r-1}}.
\end{split}
\end{equation*}
Then we define a $\U(\e)$-module
\begin{equation}\label{eq:KR module}
\W_{s,\e}^{(r)} := {\rm Im}\check{R}_r.
\end{equation}
It is proved in \cite{KO} that $\W_{s,\e}^{(r)}$ is irreducible when $\e=\e_{M|N}$, where the proof uses the crystal base of polynomial representation of $\U_{M|N}(\e)$. Now we give another proof of the irreducibility of $\W_{s,\e}^{(r)}$, which is available for arbitrary $\e$.

\begin{thm}\label{thm:W_s^{(r)}}
Let $r, s\geq 1$ be given. Then $\W_{s,\e}^{(r)}$ is non-zero if and only if $(s^r)\in \cP_{M|N}$. In this case, $\W_{s,\e}^{(r)}$ is irreducible, and it is isomorphic to $V_\e((s^r))$ as a $\ring\U(\e)$-module.


\end{thm}
\pf Let us take a sequence $\e''=(\e''_1,\dots,\e''_{n''})$ of $0,1$'s satisfying the following:
\begin{itemize}
\item[(1)] $\e$ is a subsequence of $\e''$, 

\item[(2)] we have as a $\ring\U(\e'')$-module
\begin{equation}\label{eq:decomp w.r.t. e''}
V_{\e''}((s))^{\ot r}\cong 
\bigoplus_{\la\in\cP} V_{\e''}(\la)^{\oplus K_{\la (s^r)}},
\end{equation}
where $K_{\la (s^r)}$ is the Kostka number associated to $\la$ and $(s^r)$ (cf.~Remark \ref{rem:character of poly}),

\item[(3)] if $\e'=\e_{M''|0}$ with $M''=|\{\,i\,|\,\e''_i=0\,\}|$, then we have as a $\ring\U(\e')$-module
\begin{equation}\label{eq:decomp w.r.t. e'}
V_{\e'}((s))^{\ot r}\cong 
\bigoplus_{\la\in\cP} V_{\e'}(\la)^{\oplus K_{\la (s^r)}}.
\end{equation}
\end{itemize}
Let us define a $\U(\e'')$-module $\W_{s,\e''}^{(r)}$ by the same way as in \eqref{eq:KR module}, where $\check{R}''_r$ and $V''_x$ denote the corresponding ones. We define $\W_{s,\e'}^{(r)}$, $\check{R}'_r$ and $V'_x$ similarly.

By Lemma \ref{lem:comm diagram of R matrix}, we have the following commutative diagram:
\begin{equation*}\label{eq:comm diagram of R matrix-2}
\xymatrixcolsep{4pc}\xymatrixrowsep{2pc}\xymatrix{
V''_{q^{r-1}}\otimes \cdots\otimes V''_{q^{1-r}} \ar@{->}^{\check{R}''_r}[r] \ar@{->}^{\mf{tr}^{\e''}_{\e'}}[d] &\ V''_{q^{1-r}}\otimes \cdots\otimes V''_{q^{r-1}} \ar@{->}^{\mf{tr}^{\e''}_{\e'}}[d] \\  
V'_{q^{r-1}}\otimes \cdots\otimes V'_{q^{1-r}} \ar@{->}^{\check{R}'_r}[r] &\  V'_{q^{1-r}}\otimes \cdots\otimes V'_{q^{r-1}}}\quad 
\end{equation*}

By \eqref{eq:decomp w.r.t. e''}, \eqref{eq:decomp w.r.t. e'} and Proposition  \ref{prop:truncation of poly}, the decomposition of $\W_{s,\e''}^{(r)}$ into polynomial $\ring\U(\e'')$-modules is the same as that of $\W_{s,\e'}^{(r)}$ into polynomial $\ring\U(\e')$-modules. It is well-known that $\W_{s,\e'}^{(r)}$ is irreducible and isomorphic to $V_{\e'}((s^r))$ as a $\ring\U(\e')$-module since $\U(\e'')\cong U_{q}(A_{M''-1}^{(1)})$. 
Therefore, $\W_{s,\e''}^{(r)}$ is irreducible and isomorphic to $V_{\e''}((s^r))$ as a $\ring\U(\e'')$-module. 

Again by Lemma \ref{lem:comm diagram of R matrix}, we have the following commutative diagram:
\begin{equation*}\label{eq:comm diagram of R matrix-3}
\xymatrixcolsep{4pc}\xymatrixrowsep{2pc}\xymatrix{
V''_{q^{r-1}}\otimes \cdots\otimes V''_{q^{1-r}} \ar@{->}^{\check{R}''_r}[r] \ar@{->}^{\mf{tr}^{\e''}_{\e'}}[d] &\ V''_{q^{1-r}}\otimes \cdots\otimes V''_{q^{r-1}} \ar@{->}^{\mf{tr}^{\e''}_{\e'}}[d] \\  
V_{q^{r-1}}\otimes \cdots\otimes V_{q^{1-r}} \ar@{->}^{\check{R}_r}[r] &\  V_{q^{1-r}}\otimes \cdots\otimes V_{q^{r-1}}}\quad 
\end{equation*}
Since $\mf{tr}^{\e''}_{\e'}(V_{\e''}((s^r)))$ is non-zero if and only if $(s^r)\in \cP_{M|N}$, which is equal to $V_{\e}((s^r))$ in this case, it follows that $\W_{s,\e}^{(r)}$ is non-zero if and only if $(s^r)\in \cP_{M|N}$. This implies in this case that $\W_{s,\e}^{(r)}$ is irreducible, and it is isomorphic to $V_\e((s^r))$ as a $\ring\U(\e)$-module.
\qed

The following can be proved by similar arguments.
\begin{cor}
Suppose that $(s^r)\in \cP_{M|N}$ is given. 
\begin{itemize}
\item[(1)] If $r\leq M$ and $M\geq 3$, then $\mf{tr}^\e_{\e'}\left(\W_{s,\e}^{(r)}\right)$ is the Kirillov-Reshetikhin module of type $A_{M-1}^{(1)}$ corresponding to the partition $(s^r)$, where $\e'=\e_{M|0}$.

\item[(2)] If $s\leq N$ and $N\geq 3$, then $\mf{tr}^\e_{\e'}\left(\W_{s,\e}^{(r)}\right)$ is the Kirillov-Reshetikhin module of type $A_{N-1}^{(1)}$ corresponding to the partition $(r^s)$, where $\e'=\e_{0|N}$.

\end{itemize}
\end{cor}

\begin{rem}{\rm
As in case of $\e=\e_{M|N}$ \cite{KO}, we also expect that $\W^{(r)}_{s,\e}$ has a crystal base for arbitrary $\e$ (cf.~Remark \ref{rem:crystal base of poly}).
 }
\end{rem}


One may use a similar argument as in the proof of Theorem \ref{thm:W_s^{(r)}} to prove the irreducibility of a tensor product of $\W_{l,\e}(x)$'s and its image under $R$ matrix in some special cases. 
Let $l_1,\dots,l_r\in \Z_+$ and $x_1,\dots,x_r\in \mathbb{Q}(q)$ be given and let $\e'=\e_{M|0}$.

\begin{prop}\label{prop:irreducibility-1}
If $M$ is sufficiently large and $\W_{l_1,\e'}(x_1)\ot \cdots \ot \W_{l_r,\e'}(x_r)$ is irreducible, then 
$\W_{l_1,\e}(x_1)\ot \cdots \ot \W_{l_r,\e}(x_r)$ is also irreducible.
\end{prop}
\pf Suppose that $\W_{l_1,\e}(x_1)\ot \cdots \ot \W_{l_r,\e}(x_r)$ is not irreducible and let $W$ be a proper non-trivial submodule. Since $M$ is sufficiently large, the multiplicity of $V_\e(\la)$ for $\la\in \cP$ in $\W_{l_1,\e}(x_1)\ot \cdots \ot \W_{l_r,\e}(x_r)$ is equal to that of 
$V_{\e'}(\la)$ for $\la\in \cP$ in $\W_{l_1,\e'}(x_1)\ot \cdots \ot \W_{l_r,\e'}(x_r)$ (cf.~Remark \ref{rem:character of poly}). This also holds for $W$ and $\mf{tr}^\e_{\e'}(W)$, which implies that $\mf{tr}^\e_{\e'}(W)$ is a proper non-zero subspace of $\mf{tr}^\e_{\e'}(\W_{l_1,\e}(x_1)\ot \cdots \ot \W_{l_r,\e}(x_r))$. 
Since $\mf{tr}^\e_{\e'}(W)=W\cap \mf{tr}^\e_{\e'}\left(\W_{l_1,\e}(x_1)\ot \cdots \ot \W_{l_r,\e}(x_r)\right)$, it follows that $\mf{tr}^\e_{\e'}(W)$ is a proper non-zero $\U(\e')$-submodule, which is a contradiction.
\qed

\begin{rem}{\rm
Proposition \ref{prop:irreducibility-1} together with the irreducibility of $\W_{l,\e'}\ot \W_{m,\e'}$ also implies Theorem \ref{thm:main theorem-1} when $M\geq 3$. But we do not know whether it holds for $M=2$. We also would like to point out that the proof of Theorem \ref{thm:main theorem-1} has its own interest since it describes a new connected crystal graph structure on $\mc{B}_{l,\e}\ot \mc{B}_{m,\e}/\{\pm 1\}$.
}
\end{rem}

\begin{prop}\label{prop:irreducibility-2}
Suppose that $x_i/x_{i+1}\not\in q^{-2\Z_+}$ for $1\leq i\leq r-1$.  
If $M$ is sufficiently large and the image of 
$$\check{R}'_{w_0}(x_1,\dots,x_r) : \W_{l_1,\e'}(x_1)\ot \cdots \ot \W_{l_r,\e'}(x_r) \longrightarrow \W_{l_r,\e'}(x_r)\ot \cdots \ot \W_{l_1,\e'}(x_1)$$ is irreducible, then 
the image of 
$$\check{R}_{w_0}(x_1,\dots,x_r) : \W_{l_1,\e}(x_1)\ot \cdots \ot \W_{l_r,\e}(x_r) \longrightarrow \W_{l_r,\e}(x_r)\ot \cdots \ot \W_{l_1,\e}(x_1)$$ is also irreducible, where $\check{R}'_{w_0}(x_1,\dots,x_r)$ is the restriction of $\check{R}_{w_0}(x_1,\dots,x_r)$ on $\W_{l_1,\e'}(x_1)\ot \cdots \ot \W_{l_r,\e'}(x_r)$.
\end{prop}
\pf It follows from Lemma \ref{lem:comm diagram of R matrix} and the same argument as in Proposition \ref{prop:irreducibility-1}.
\qed

{\small

\end{document}